\documentclass[11pt]{amsart}
\usepackage{amsmath, amsthm, amssymb}
\usepackage{graphicx, array}
\usepackage{longtable}
\usepackage{epstopdf}

\newtheorem{Thm}{Theorem}[section]
\newtheorem{Prop}[Thm]{Proposition}
\newtheorem{Lem}[Thm]{Lemma}

\theoremstyle{definition}

\newtheorem{Rem}[Thm]{Remark}
\newtheorem{Ex}[Thm]{Example}

\newenvironment{Proof}{\rm \trivlist\item[\hskip \labelsep{\bf
Proof.\quad}]}{\hfill\qed\par\medskip\endtrivlist}

\newcommand{\X}{X \cup X^{-1}} 
\newcommand{\comp}{\mathcal{C}} 
\newcommand{\compb}{\mathcal{D}} 
\newcommand{\graph}{\Gamma} 
\newcommand{\graphb}{\Delta} 
\def\M#1{M_{#1}} 
\def\init#1{\alpha(#1)} 
\def\ter#1{\omega(#1)} 
\def\pcat#1{\operatorname{IC}(#1)} 
\def\cat{C} 

\newcommand{\Stab}{\operatorname{Stab}}
\newcommand{\FIM}{\operatorname{FIM}}
\newcommand{\SIM}{\operatorname{SIM}}
\newcommand{\bd}{\partial}
\newcommand{\bl}{bl}

\def\bound#1{bw(#1)}

\title{Inverse monoids and immersions of 2-complexes}

\author{John Meakin}
\email{jmeakin@math.unl.edu}
\address{Department of Mathematics\\
University of Nebraska-Lincoln\\
Lincoln, NE 68588, USA}

\author{N\'ora Szak\'acs }
\email{szakacsn@math.u-szeged.hu}
\address{Bolyai Institute\\
University of Szeged\\
Aradi v\'ertan\'uk tere 1.\\
H-6720 Szeged, Hungary}
\thanks{This research was realized in the frames of T\'AMOP 4.2.4.
A/2-11-1-2012-0001 ``National Excellence Program - Elaborating and
operating an inland student and researcher personal support system
convergence program". The project was subsidized by the European
Union and co-financed by the European Social Fund.
This research was partially supported by the Hungarian National Foundation
for Scientific Research grant no. K104251.}

\begin{document}

\begin{abstract}
It is well known that under mild conditions on a connected
topological space $\mathcal X$, connected covers of $\mathcal X$ may
be classified via conjugacy classes of subgroups of the fundamental
group of $\mathcal X$. In this paper, we extend these results to the
study of {\em immersions} into $2$-dimensional $CW$-complexes. An
immersion $f : {\mathcal D} \rightarrow \mathcal C$ between
$CW$-complexes is a cellular map such that each point $y \in
{\mathcal D}$ has a neighborhood $U$ that is mapped homeomorphically
onto $f(U)$ by $f$. In order to classify immersions into a
$2$-dimensional $CW$-complex $\mathcal C$, we need to replace the
fundamental group of $\mathcal C$ by an appropriate inverse monoid.
We show how conjugacy classes of the closed inverse submonoids of
this inverse monoid may be used to classify connected immersions
into the complex.

\end{abstract}

\maketitle

Dedicated to Stuart Margolis, on the occasion of his 60th birthday.

\section{Introduction}

It is well known that under mild restrictions on a topological space
$\mathcal X$, connected covers of $\mathcal X$ may be classified via
conjugacy classes of subgroups of the fundamental group of $\mathcal
X$. For this fact, and for general background in topology, we refer
to the book by Munkres \cite{munk}.

In this paper we study connected immersions between
finite-dimensional $CW$-complexes. A $CW$-complex $\comp$  is
obtained from a discrete set ${\comp}^0$ (the $0$-skeleton of
$\comp$) by iteratively attaching cells of dimension $n$ to the
$(n-1)$-skeleton ${\comp}^{n-1}$ of $\comp$  for $n \geq 1$. We
refer the reader to Hatcher's text \cite{hatch}, for the precise
definition and basic properties of $CW$-complexes. In particular a
continuous map between $CW$-complexes is homotopic to a {\em
cellular map} (\cite{hatch}, Theorem 4.8),  that is a continuous
function that maps cells  to cells of the same or lower dimension,
so we will regard maps between $CW$-complexes as cellular maps. A
subcomplex of a $CW$-complex is a closed subspace that is a union of
cells.

An {\em immersion} of a $CW$-complex $\mathcal D$  into a
$CW$-complex $\comp$ is a cellular map $f : {\mathcal D} \rightarrow
{\comp}$ such that each point $y \in \mathcal D$ has a neighborhood
$U$ which is mapped homeomorphically onto $f(U)$ by $f$. So $f$ maps
$n$-cells to $n$-cells. Thus if $\comp$ is an $n$-dimensional
$CW$-complex, then $\mathcal D$ is an $m$-dimensional $CW$-complex
with $m \leq n$. Every subcomplex of an $n$-dimensional $CW$-complex
$\comp$ immerses into $\comp$. Every covering space of a
$CW$-complex $\comp$ has a $CW$-complex structure, and every
covering map is in particular an immersion.

We classify connected immersions into a $2$-dimensional $CW$-complex
$\comp$ via conjugacy classes of closed inverse submonoids of a
certain inverse monoid associated with $\comp$.  The closed inverse
submonoids of this inverse monoid enable us to keep track of the
$1$-cells and $2$-cells of $\comp$ that lift under the immersion, in
much the same way as the subgroups of the fundamental group of
$\comp$ enable us to encode coverings of $\comp$. We provide an
iterative process for constructing the immersion associated with a
closed inverse submonoid of this inverse monoid. In many cases this
iterative procedure provides an algorithm for constructing the
immersion, in particular if the closed inverse submonoid is finitely
generated and $\comp$ has finitely many $2$-cells.

Section 2 of the paper outlines basic  material on presentations of
inverse monoids  that we will need to build an inverse monoid
associated with a  $2$-complex $\comp$. Section 3 describes an
iterative procedure for constructing closed inverse submonoids of an
inverse monoid from generators for the submonoid. Section 4
describes connected immersions between $2$-complexes in terms of a
labeling of $1$-cells and $2$-cells. The main results of the paper
linking immersions over a $2$-complex $\comp$ and closed inverse
submonoids of an inverse monoid associated with $\comp$ are
described  in detail in Section 5 of the paper (Theorem
\ref{premain}, Theorem \ref{main} and Theorem \ref{mainalg}). We
close in Section 6 with several examples illustrating the
connections between immersions over $2$-complexes and the associated
closed inverse submonoids.

These results extend some work of Margolis and Meakin \cite{MM1}
that classifies connected immersions over graphs ($1$-dimensional
$CW$-complexes) via closed inverse submonoids of free inverse
monoids.  Some related work may be found in the thesis of Williamson
\cite{will}. However,  the notion of immersion in this paper is
considerably more general than the notion of immersion between
$2$-complexes in  \cite{will}.

\section{$X$-graphs and inverse monoids}

Let $X$ be a  set and $X^{-1}$ a disjoint set in one-one
correspondence with $X$ via a map $x \rightarrow x^{-1}$ and define
$(x^{-1})^{-1} = x$. We extend this to a map on $(\X)^*$ by defining
$(x_{1}x_{2} \cdots x_{n})^{-1} = x_{n}^{-1} \cdots
x_{2}^{-1}x_{1}^{-1}$, giving $(\X)^*$ the structure of the free
monoid with involution on $X$. Throughout this paper by an {\em
$X$-graph} (or just an {\em edge-labeled graph} if the labeling set
$X$  is understood) we mean a strongly connected digraph $\graph$
with edges labeled over the set $\X$ such that the labeling is
consistent with an involution: that is, there is an edge labeled $x
\in \X$ from vertex $v_1$ to vertex $v_2$ if and only if there is an
inverse edge labeled $x^{-1}$ from $v_2$ to $v_1$. The initial
vertex of an edge $e$ will be denoted by $\init e$ and the terminal
vertex by $\ter e$. If $X = \emptyset$, then we view $\graph$ as the
graph with one vertex and no edges.

The label on an edge $e$ is denoted by $\ell(e) \in \X$. There is an
evident notion of {\em path} in an $X$-graph. A path $p$ with
initial vertex $v_1$ and terminal vertex $v_2$ will be called a
$(v_{1},v_{2})$ path. The initial (resp. terminal) vertex of a path
$p$ will be denoted by ${\alpha}(p)$ (resp. ${\omega}(p)$). The
label on the path $p = e_{1}e_{2} \ldots e_k$ is the word $\ell(p) =
\ell(e_{1})\ell(e_{2})\ldots \ell(e_{k}) \in (\X)^*$.

It is customary when sketching diagrams of such graphs to include
just the positively labeled edges (with labels from $X$) in the
diagram.

$X$-graphs occur frequently in the literature.  For example, the
bouquet of $|X|$ circles is the $X$-graph  $B_{X}$ with one vertex
and one positively labeled edge labeled by $x$ for each $x \in X$.
The Cayley graph ${\graph}(G,X)$ of a group $G$ relative to a set
$X$ of generators is an $X$-graph: its vertices are the elements of
$G$ and it has an edge labeled by $x$ from $g$ to $gx$ for each $x
\in \X$.

 If we designate an initial vertex (state) $\alpha$ and a
terminal vertex (state) $\beta$ of $\graph$, then the birooted
$X$-graph ${\mathcal A} = ({\alpha}, {\graph}, {\beta})$ may be
viewed as an automaton. See for example the book of Hopcroft and
Ullman \cite {HU} for basic information about automata theory. The
language accepted by this automaton is the subset $L({\mathcal A})$
of $(\X)^*$ consisting of the words in $(\X)^*$ that label paths in
$\graph$ starting at $\alpha$ and ending at $\beta$. This automaton
is called an {\em inverse automaton} if it is deterministic (and
hence co-deterministic), i.e. if for each vertex $v$ of $\graph$
there is at most one edge with a given label starting or ending at
$v$. A deterministic $X$-graph $\graph$ determines an immersion of
$\graph$ into $B_X$, obtained by mapping an edge labeled by $x \in
\X$ onto the corresponding edge in $B_X$.

Recall that an {\em inverse monoid} is a monoid $M$ with the
property that for each $a \in M$ there exists a unique element
$a^{-1} \in M$ (the inverse of $a$) such that $a =aa^{-1}a$ and
$a^{-1} = a^{-1}aa^{-1}$.  Every inverse monoid may be embedded in a
suitable symmetric inverse monoid $\SIM(V)$. Here $\SIM(V)$ is the
monoid of all partial injective functions from  $V$ to $V$ (i.e.
bijections between subsets of the set $V$) with respect to the usual
composition of partial maps. If $\graph$ is a deterministic
$X$-graph, then each letter $x \in \X$ determines a partial
injection of the set $V$ of vertices of $\graph$ that maps a vertex
$v_1$ to a vertex $v_2$ if there is an edge labeled by $x$ from
$v_1$ to $v_2$. The submonoid of $\SIM(V)$ generated by these
partial maps is an inverse monoid, called the transition monoid of
the graph $\graph$.

We refer the reader to the books by Lawson \cite{law} or Petrich
\cite{Pet} for the basic theory of inverse monoids. In particular,
the {\em natural partial order} on an inverse monoid $M$ is defined
by $a \leq b$ iff $a = eb$ for some idempotent $e \in M$, or equivalently, if $a=aa^{-1}b$. This
corresponds to restriction of partial injective maps when $M =
\SIM(V)$. Note that if $M$ is a group, than the natural partial order is just the equality. See \cite{law} or \cite{Pet} for the important role that the natural
partial order plays in the structure of inverse monoids.
If we factor an inverse monoid $M$ by the congruence generated by pairs of the form $(aa^{-1},1)$,
 $a \in M$, we obtain a group. This congruence is denoted by $\sigma$, and $M/\sigma$ is
 in fact the \emph{greatest group homomorphic image of} $M$.

Since inverse monoids form a variety of algebras (in the sense of
universal algebra - i.e. an equationally defined class of algebras),
free inverse monoids exist. We will denote the free inverse monoid
on a set $X$ by $\FIM(X)$. This is the quotient of $(\X)^*$, the
free monoid with involution, by the congruence that identifies
$ww^{-1}w$ with $w$ and $ww^{-1}uu^{-1}$ with $uu^{-1}ww^{-1}$ for
all words $u,w \in (\X)^*$. See \cite{Pet} or \cite{law} for much
information about $\FIM(X)$. In particular, \cite{Pet} and
\cite{law} provide an exposition of  Munn's solution \cite{munn} to
the word problem for $\FIM(X)$ via birooted edge-labeled trees
called {\em Munn trees}.

In his thesis \cite{Ste2} and paper \cite{Ste1}, Stephen initiated
the theory of presentations of inverse monoids by extending Munn's
results about free inverse monoids to arbitrary presentations of
inverse monoids. Here, a presentation of an inverse monoid $M$,
denoted $M = Inv \langle X\ |\ u_{i} = v_{i}, i \in I \rangle$ (where
the $u_i$ and $v_i$ are words in $(\X)^*$) is the quotient of
$\FIM(X)$ obtained by  imposing the relations $u_{i} = v_{i}$ in the
usual way. In order to study the word problem for such
presentations, Stephen considers the {\em Sch\"utzenberger graph}
$S{\graph}(M,X,w)$ (or simply $S{\graph}(w)$ if the presentation is
understood) of each word $w \in (\X)^*$. The Sch\"utzenberger graph of $w$ is the
restriction of the Cayley graph of $M$ to the $\mathcal R$-class of
$w$ in $M$. That is, the vertices of $S{\graph}(w)$ are the elements
$u \in M$ such that $uu^{-1} = ww^{-1}$ in $M$; there is an edge
labeled by $x \in \X$ from $u$ to $v$ if $uu^{-1} = vv^{-1} =
ww^{-1}$ and $ux = v$ in $M$. (Here, for simplicity of notation, we
are using the same notation for a word $w \in (\X)^*$ and its
natural image in $M$; the context guarantees that no confusion
should occur.)

The Sch\"utzenberger graphs of $M$ are just the strongly connected
components of the Cayley graph of $M$ relative to the set $X$ of
generators for $M$. Of course, if $G$ is a group, then it has just
one Sch\"utzenberger graph, which is the Cayley graph
${\graph}(G,X)$.
The {\em Sch\"utzenberger automaton} $S{\mathcal A}(w)$ of a word $w
\in (\X)^*$ is defined to be the birooted $X$-graph $S{\mathcal
A}(w) = (ww^{-1},S{\graph}(w),w)$. Thus $S{\mathcal A}(w)$ is an
inverse automaton. In his paper \cite{Ste1}, Stephen proves the
following result.

\begin{Thm}

Let $M = Inv \langle X \ |\ u_{i} = v_{i}, i \in I \rangle$ be a
presentation of an inverse monoid. Then

(a) For each word $u \in (\X)^*$, the  language accepted by the \\
Sch\"utzenberger automaton $S{\mathcal A}(u)$ is the set of all
words $w \in (\X)^*$ such that $u \leq w$ in the natural partial
order on $M$.

(b) $u = w$ in $M$ iff $u \in L(S{\mathcal A}(w))$ and $w \in
L(S{\mathcal A}(u))$.

(c) The word problem for $M$ is decidable iff there is an algorithm
for deciding membership in $L(S{\mathcal A}(w))$ for each word $w
\in (\X)^*$.

\end{Thm}

\section{Closed inverse submonoids of inverse monoids}

For each subset $N$ of an inverse monoid $M$, we denote by
$N^{\omega}$ the set of all elements $m \in M$ such that $m \geq n$
for some $n \in N$.  The subset $N$ of $M$ is called {\em closed} if
$N = N^{\omega}$. Thus the image in $M$ of the language accepted by
a Sch\"utzenberger automaton $S{\mathcal A}(u)$ of a word $u$
relative to a presentation of $M$ is a closed
subset of $M$.

Closed inverse submonoids of an inverse monoid $M$ arise naturally
in the representation theory of $M$ by partial injections on a set
 \cite{schein}. An inverse monoid $M$ acts (on the right) by
 injective partial functions on a set $Q$ if there is a homomorphism
 from $M$ to $\SIM(Q)$. Denote by $qm$ the image of $q$ under the
 action of $m$ if $q$ is in the domain of the action by $m$.
 The following basic fact is well known (see \cite{schein}).

 \begin{Prop}
\label{stabclosed}

 If $M$ acts on $Q$ by injective partial functions, then for every $q \in
 Q, \, Stab(q) = \{m \in M : qm = q\}$ is a closed inverse submonoid
 of $M$.

 \end{Prop}

 Conversely, given a  closed inverse submonoid $H$ of $M$, we
 can construct a transitive representation of $M$ as follows. A
 subset of $M$ of the form $(Hm)^{\omega}$ where $mm^{-1} \in H$ is
 called a {\em right ${\omega}$-coset} of $H$. Let $X_H$ denote the
 set of right $\omega$-cosets of $H$. If $m \in M$, define an action
 on $X_H$ by $Y \cdot m = (Ym)^{\omega}$ if $(Ym)^{\omega} \in X_H$ and
 undefined otherwise. This defines a transitive action of $M$ on
 $X_H$. Conversely, if $M$ acts transitively on $Q$, then this action
 is equivalent in the obvious sense to the action of $M$ on the
 right $\omega$-cosets of $Stab(q)$ in $M$ for any $q \in Q$. See
 \cite{schein} or \cite{Pet} for details.

 The {\em $\omega$-coset graph} ${\graph}_{(H,X)}$ (or just ${\graph}_H$ if $X$ is understood)
 of a closed inverse
 submonoid $H$ of an $X$-generated inverse monoid $M$  is
 constructed as follows. The
 set of vertices of ${\graph}_H$ is $X_H$ and there is an edge
 labeled by $x \in \X$ from $(Ha)^{\omega}$ to $(Hb)^{\omega}$ if
 $(Hb)^{\omega} = (Hax)^{\omega}$.
 Then ${\graph}_H$ is a
 deterministic
 $X$-graph. The
 birooted $X$-graph $(H,{\graph}_{H},H)$ is called the {\em
 $\omega$-coset automaton} of $H$. The language accepted by this
 automaton is $H$ (or more precisely the set of words $w \in (\X)^*$
 whose natural image in $M$ is in $H$).
 Clearly, if $G$ is a group generated by $X$, then ${\graph}_H$ coincides with the coset graph of
 the
 subgroup $H$ of $G$.

 Let $M$ be an inverse
monoid given by a presentation $M = Inv\langle X\ |\ u_{i}=v_{i}, \, i
\in I\rangle$, and let $Y$ be a subset of $(X \cup X^{-1})^*$. Let
$\langle Y \rangle ^{\omega} $ denote the closed inverse submonoid
of $M$ generated by  the natural image of $Y$ in $M$.
We now provide an iterative construction of the $\omega$-coset
automaton of $\langle Y \rangle ^{\omega} $. The construction
extends the well-known construction of Stallings \cite{stall}  of a
finite graph
 associated with each finitely generated subgroup of a free group.
See also \cite{MM1} for the automata-theoretic point of view on
Stallings' construction.

In \cite{Ste2}, Stephen shows that the class of all birooted
$X$-graphs  forms a cocomplete category, and hence directed systems
of birooted $X$-graphs have direct limits in this category. See Mac
Lane \cite{McL} for background in category theory.  Morphisms in
this category are graph morphisms that take edges to edges and
preserve edge labelings and initial (terminal) roots.

Given a finite presentation $M = Inv \langle X\ |\ u_{i} = v_{i},\ i =
1, \ldots ,n \rangle$ of an inverse monoid, we consider two types of
operations on $X$-graphs (or birooted $X$-graphs), namely edge
foldings (in the sense of Stallings \cite{stall}) and expansions. If
$e_1$ and $e_2$ are two edges with the same label and the same initial
or terminal vertex, then an edge folding identifies these edges (an
edge folding is called a ``determination" in Stephen's terminology
\cite{Ste1, Ste2}). Clearly, each edge folding of an $X$-graph
results in another $X$-graph. If $\graph$ is an $X$-graph with two
vertices $a$ and $b$ and a path from $a$ to $b$ labeled by one side
(say $u_{i}$) of one of the defining relations $u_{i} = v_{i}$ of
the monoid $M$, but no path labeled by the other side, then we expand
$\graph$ to create another $X$-graph $\graphb$ by adding a new path
from $a$ to $b$ labeled by the other side ($v_i$) of the relation.
One of the results of Stephen \cite{Ste1} (Lemma 4.7) is that these
processes are confluent.

The set of birooted $X$-graphs obtained by applying successive
expansions and edge foldings to a birooted $X$-graph ${\mathcal A} =
({\alpha}, {\graph}, {\beta})$ forms a directed system in the
category of birooted $X$-graphs. The direct limit (colimit) of this
system is an inverse automaton that we will denote by ${\mathcal
A}^{\omega}$. This automaton is {\em complete}, in the sense that no
edge foldings or expansions may be applied. Of course if finitely
many applications of edge foldings and expansions transform
$\mathcal A$ into a complete automaton $\mathcal B$, then ${\mathcal
B} = {\mathcal A}^{\omega}$.

Any automaton $\mathcal A'$ obtained from $\mathcal A$ by applying
successive expansions and edge foldings is called an {\em
approximate automaton} of ${\mathcal A}^{\omega}$.

\begin{Thm}

Let $M = Inv \langle X : u_{i} = v_{i}, \, i = 1, \ldots n \rangle$
be a finitely presented inverse monoid. If ${\mathcal A}$ is a
birooted $X$-graph (i.e. automaton) accepting the language $L
\subseteq (\X)^*$, then the language accepted by the direct limit
automaton ${\mathcal A}^{\omega}$ is $L^{\omega} = \{w \in (X \cup
X^{-1})^{*} : w \geq s$ in $M$ for some $s \in L\}$.

\end{Thm}

\begin{Proof} The proof follows by a modification of the proof of
Theorem 4.12 of Stephen \cite{Ste2}, where it is proved that the
Sch\"utzenberger automaton $S{\mathcal A}(s)$ of a word $s \in
(\X)^*$ is the colimit $Lin(s)^{\omega}$, where $Lin(s)$ is the
``linear automaton" of $s$. See also Theorem 5.10 of \cite{Ste1} for
a closely related result. The basic idea of the proof is that
application of an expansion to some automaton $\mathcal A'$ just
augments the language $L({\mathcal A'})$ by words that are equal in
$M$ to words in $L({\mathcal A'})$, while an edge folding augments
this language by words that are greater than or equal in $M$ to
words in $L({\mathcal A'})$. We provide some more detail below.

Let ${\mathcal A}^{\omega} = ({\alpha}^{\omega}, {\graph}^{\omega},
{\beta}^{\omega})$. If $w \in L({\mathcal A}^{\omega})$, then the
path labeled by $w$ lifts to a path labeled by $w$ from ${\alpha}'$
to ${\beta}'$ in some approximate automaton $\mathcal A' =
({\alpha}',{\graph}',{\beta}')$ of ${\mathcal A}^{\omega}$ by
Theorem 2.11 of \cite{Ste2}.  This implies that  $w \in L(\mathcal
A')$. But it follows as in the proof of Theorem 5.5 and Lemma 5.6 of
\cite{Ste1} that if $\mathcal A'$ is an approximate automaton of
${\mathcal A}^{\omega}$, then $L \subseteq L({\mathcal A'})
\subseteq L^{\omega}$. Hence $L({\mathcal A}^{\omega}) \subseteq
L^{\omega}$.

Conversely, if $w \geq s$ for some $s \in L$, then by Theorem 2.1
above, $w \in L(S{\mathcal A}(s))$. So $w$ is in the language
accepted by some approximate automaton $\mathcal B'$ of $S{\mathcal
A}(s)$ by Theorem 5.12 of \cite{Ste1}. The automaton $\mathcal B'$
is obtained from the linear automaton of $s$ by a finite number of
edge foldings and expansions. Since $s \in L = L({\mathcal A})$, we
may apply the same sequence of edge foldings and expansions to
$\mathcal A$ to obtain an approximate automaton $\mathcal A'$ of
${\mathcal A}^{\omega}$, and hence $w$ is in the language accepted by
this approximate automaton $\mathcal A'$.
Since there is a
morphism from $\mathcal A'$ to ${\mathcal A}^{\omega}$ by definition
of the colimit, it follows from Lemma 2.4 of \cite{Ste1} that $w \in
L({\mathcal A}^{\omega})$.

\end{Proof}

We now apply Stephen's iterative process as described above to
construct the closed inverse submonoid of $M$ generated by a
 subset $Y$ of $(X \cup X^{-1})^*$. Start with the ``flower
automaton" ${\mathcal F}(Y)$. This is the birooted $X$-graph with
one distinguished state $1$ designated as initial and terminal state
and a closed path based at $1$ labeled by the word $y$ for each $y \in Y$.
(This is a finite automaton if $Y$ is finite of course.) Now
successively apply edge foldings and expansions to ${\mathcal F}(Y)$
to obtain the limit automaton ${\mathcal F}(Y)^{\omega}$.

\begin{Thm}
\label{cosetalg}

Let $M = Inv\langle X \ |\ u_{i} = v_{i},\ i = 1, \ldots , n\rangle$ be
a finitely presented inverse monoid, let $Y$ be a  subset of $(X
\cup X^{-1})^*$, and construct the inverse automaton ${\mathcal
F}(Y)^{\omega}$ obtained from the flower automaton ${\mathcal F}(Y)$
by iteratively applying the processes of edge foldings and
expansions as described above.  Then the language $L({\mathcal
F}(Y)^{\omega})$ accepted by this automaton is $\{w \in (X \cup
X^{-1})^{*} : w \in \langle Y \rangle ^{\omega}\}$, and
 ${\mathcal F}(Y)^{\omega}$ is the $\omega$-coset automaton of the closed inverse
 submonoid $\langle
Y\rangle ^{\omega} $ of $M$. Thus the membership problem for the
closed inverse submonoid $\langle Y\rangle ^{\omega} $  is decidable
if and only if there is an algorithm for deciding membership in the
language $L({\mathcal F}(Y)^{\omega})$.

\end{Thm}

\begin{Proof}

The fact that $L({\mathcal F}(Y)^{\omega}) = \{w \in (X \cup
X^{-1})^{*} : w \in \langle Y \rangle ^{\omega}\}$ is immediate from
Theorem 3.2 above. Hence the automaton ${\mathcal F}(Y)^{\omega}$
and the $\omega$-coset automaton of the closed inverse
 submonoid $\langle
Y\rangle ^{\omega} $ are birooted deterministic  $X$-graphs that
accept the same language. But  any two birooted (connected)
deterministic $X$-graphs that accept the same language are
isomorphic as birooted $X$-graphs since inverse automata are minimal
(see \cite{reut}, Lemma 1).

\end{Proof}

This theorem shows in particular that the membership problem for the
finitely generated closed inverse submonoid $\langle Y \rangle
^{\omega} $ of $M$ is decidable if the iterative procedure described
above for constructing ${\mathcal F}(Y)^{\omega}$ terminates after a
finite number of edge foldings and expansions, since in that case
${\mathcal F}(Y)^{\omega}$ is a finite inverse automaton.

We remark that if $M$ is the free group $FG(X)$, viewed as an
inverse monoid with presentation $FG(X) = Inv\langle X\ |\ xx^{-1} =
x^{-1}x = 1 \rangle $, then finitely generated closed inverse
submonoids of $M$ coincide with finitely generated  subgroups of
$FG(X)$, and the construction of ${\mathcal F}(Y)^{\omega}$ from a
finite set $Y$ of words produces the coset graph of the subgroup.
The core of this graph is, of course, the Stallings graph (automaton)
of the corresponding subgroup \cite{stall}, obtained by pruning all
trees off the coset graph;  the reduced words accepted by the coset
automaton (or by the Stallings automaton) coincide with the reduced
words in the subgroup.

\section{Immersions of 2-complexes}

Recall the following definition \cite{hatch} of a finite dimensional
$CW$-complex $\comp$:
\begin{enumerate}
    \item Start with a discrete set $\comp^0$ , the $0$-cells of $\comp$.
    \item Inductively, form the $n$-skeleton $\comp^n$ from $\comp^{n-1}$ by attaching $n$-cells
    $C^n_\alpha$ via maps
$\varphi_\alpha \colon S^{n-1} \to \comp^{n-1}$. This means that
$\comp^n$ is the quotient space of $\comp^{n-1}\ \dot\cup_\alpha\
D^n_\alpha$ under the identifications $x \sim \varphi_\alpha (x)$
for $x \in \bd D_\alpha^n$. The cell $C^n_\alpha$ is a homeomorphic
image of $D^n_\alpha - \bd D^n_\alpha$ under the quotient map.
     \item Stop the inductive process after a finite number of steps to obtain
a finite dimensional $CW$-complex $\comp$.
 \end{enumerate}
 The dimension of the
complex is the largest dimension of one of its cells. We denote the
set of $n$-cells of $\comp$ by $\comp^{(n)}$. Throughout the
remainder of this paper, by a {\em $2$-complex} we mean a connected
$CW$-complex of dimension less than or equal to $2$. The
$1$-skeleton of a $2$-complex is an undirected graph, but it is more
convenient for our purposes to regard it as a digraph, with two
oppositely directed edges for each undirected edge.

An immersion between $CW$-complexes always maps $n$-cells to
$n$-cells, and the restriction of an immersion to a subcomplex is
also an immersion. It is easy to see that a cellular map $f \colon
\comp \to \compb$ is an immersion if and only if it is locally
injective at the $0$-cells, that is, each $0$-cell $v \in
\comp^{(0)}$ has a neighborhood that is homeomorphic to its image
under $f$. For graphs, this definition of immersions is equivalent to
Stallings' definition in \cite{stall}.

In this section, we classify immersions over $2$-complexes using
inverse monoids. Our results extend the results of \cite{MM1}, where
the authors classify immersions over graphs by keeping track of
which closed paths lift to closed paths. This is essentially what we
do in  this paper, with the added information about when $2$-cells
lift. It will be convenient to label the $1$-cells over some set
$\X$ and the $2$-cells over some disjoint set $P$ as described
below.
With every $2$-cell, we associate a distinguished vertex (root) and
walk on its boundary, consistent with the labeling. We first
describe the process of choosing a root and boundary walk for
$2$-cells.

Let $\comp$ be a  $2$-complex and let $C$ be a $2$-cell of $\comp$
with the attaching map $\varphi_C : S^{1} \rightarrow {\comp}^1$.
Choose a point $x_0$ on the circle $S^1$ in such a way that
$\varphi_C$ maps $x_0$ to a $0$-cell of $\comp$. Consider
 $S^1$ as the interval $[0,1]$ with its endpoints
glued together and identified with $x_0$. Consider the closed path
(in the topological sense) $p_C \colon [0,1] \to {\comp}$, with $p_C
|_{(0,1)} = \varphi_C |_{(0,1)}$, $p_C(0)=p_C(1)=\varphi_C(x_{0})$.
Since the closure of every $2$-cell meets only finitely many
$0$-cells or $1$-cells (\cite{hatch}, Proposition
A.1), the image of this path corresponds to a
closed path in ${\comp}^1$ (in the graph theoretic sense) that we
call the boundary walk of $C$: we denote it by $\bound C$.
We allow for the possibility that $\bound C$ might have no edges. We
call the $0$-cell $\varphi_C(x_{0})$ the {\em base} or \emph{root}
of the $2$-cell $C$ and of the closed path $\bound C$ and denote it
by $\init C$.

Let $B_X$ be the bouquet of $|X|$ circles. We build a $2$-complex
$B_{X,P}$ by attaching labeled  $2$-cells to $B_X$ with labels
coming from a set $P$ (which we assume to be disjoint from $\X$),
and with a specified  boundary walk for each $2$-cell, as described
above. The labeling is chosen so that different $2$-cells in
$B_{X,P}$ have different labels (even if they have the same boundary
in $B_X$). We allow for the possibility that $P = \emptyset$ or that
$X = \emptyset$. Denote the label of a $2$-cell $C$ in $B_{X,P}$ by
$\ell (C) \in P$.

Every $2$-complex $\comp$ admits an immersion $f \colon \comp \to
B_{X,P}$ for some sets $X$ and $P$: one could choose $X$ as an index
set for the (undirected) edges of $\comp$ and $P$ as an index set
for the $2$-cells
for example, but
we would normally choose smaller sets $X$ and $P$ if possible. This
mapping $f$  induces a labeling on $\comp$ by giving each $1$-cell
or $2$-cell in $\comp$ the label of its image in $B_{X,P}$ under
$f$. From now on, by a {\em labeled $2$-complex}, we mean a labeling
induced by an immersion into some complex $B_{X,P}$. The
$1$-skeleton of a $2$-complex $\comp$ labeled this way is a
deterministic $X$-graph that immerses via the restriction of $f$
into $B_X$; $2$-cells of $\comp$ have the same label in $P$ if they
map to the same $2$-cell
in
$B_{X,P}$.

\begin{Ex}
\label{torus} Let $X=\{a,b\}$, $P=\{\rho\}$, and let $B_{X,P}$ be
the $2$-complex with one $2$-cell $C$ (labeled by $\rho$)
corresponding to the attaching map that takes $S^1$ to the closed
path labeled by $aba^{-1}b^{-1}$. Then $\ell (\bound
C)=aba^{-1}b^{-1}$, and $B_{X,P}$ is the presentation complex of the
free abelian group of rank $2$, and is homeomorphic to the torus. We
could have chosen any cyclic conjugate of $aba^{-1}b^{-1}$ or its
inverse and obtained the same $2$-complex, but with a different
boundary walk.

\begin{center}
\includegraphics[width=0.5\linewidth]{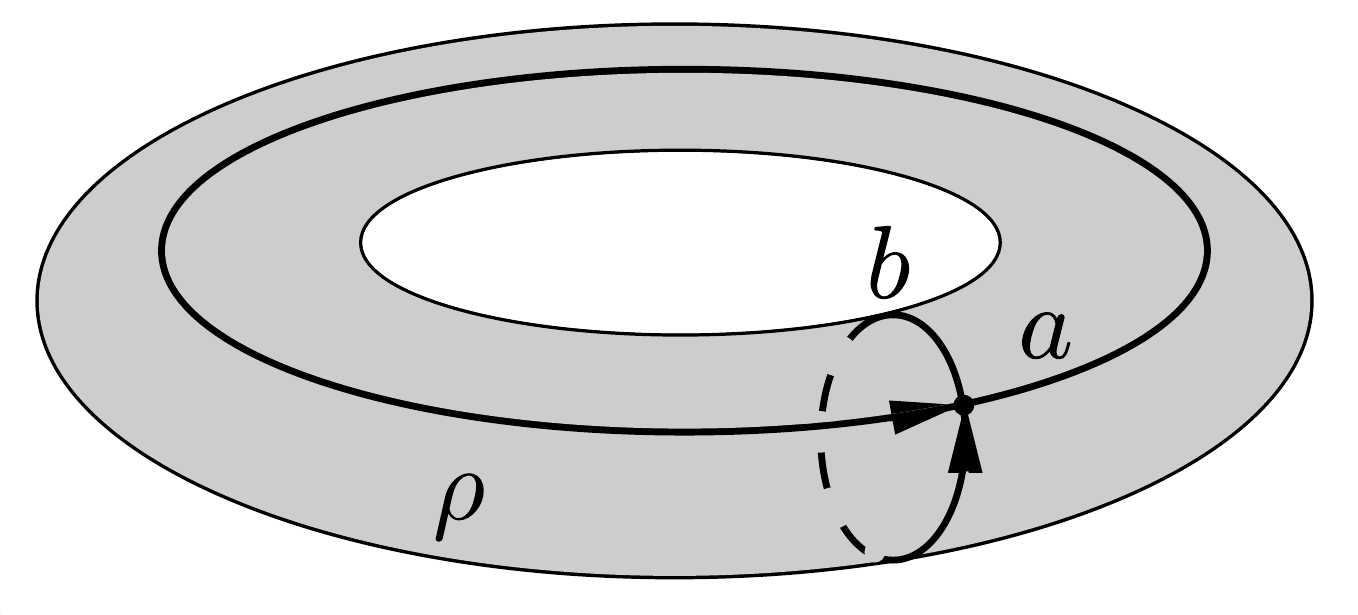}
\end{center}

\end{Ex}

If $\comp, \compb$ are $2$-complexes and $f \colon \comp \to \compb$
 an immersion, and $\compb$ is labeled by an immersion $g \colon
\compb \to B_{X,P}$, then $g \circ f \colon \comp \to B_{X,P}$ is an
immersion, and it induces a labeling on $\comp$ that is respected by
$f$; that is, $\ell (C) = \ell (f(C))$ and $\ell (e) = \ell (f(e))$
for all $2$-cells $C$ and $1$-cells $e$ in $\comp$.
 Therefore we may, without loss of
generality, assume that immersions respect the labeling.

\begin{Lem}
\label{attaching} Let $\comp, \compb$ be labeled $2$-complexes and
let $f \colon \comp \to \compb$ be an immersion that respects the
labeling. For an arbitrary $2$-cell $C$ of $\comp$,
$f(\alpha(C))=\alpha(f(C))$ and $f(\bound C)=\bound{f(C)}$.
Furthermore, $\bound C$ is uniquely determined by $f$ and $\bound
{f(C)}$.
\end{Lem}

\begin{Proof}
Let $\varphi_C \colon S^1 \to \comp$ and
$\varphi_{f(C)} \colon S^1 \to \compb$ be the attaching maps
corresponding to $C$ and $f(C)$. If $\varphi_{f(C)}$ maps the circle to a point, then so
does $\varphi_C$, and our statement
trivially holds. For the remainder of the proof, we suppose that is not the case.

We first prove that $f \circ \varphi_C=\varphi_{f(C)}$.
Consider $\comp$ as $\comp^1\ \dot{\cup}\ D_\alpha^2$ with
identifications $x \sim \varphi_\alpha (x)$ for $x \in
\bd D_\alpha^2$. Thus the closure of our $2$-cell $C$ is a quotient
of $\varphi_C(S^1)\ \dot\cup\ D^2$ by identifying the points $x
\sim \varphi_C (x)$ for $x \in \bd D^2$. Since $f$ is an
immersion, then $f|_C$ is a homeomorphism, and so $f(\overline C)$
is $f(\varphi_C(S^1))\ \dot\cup\ D^2$ with the identifications $x
\sim f(\varphi_C (x))$ for $x \in \bd D^2 $. But $f(\overline C)$
is the closure of the $2$-cell $f(C)$ in $\compb$, so it is also
$\varphi_{f(C)}(S^1)\ \dot\cup\ D^2$ with identifications $x \sim
\varphi_{f(C)} (x)$ for $x \in \bd D^2$. That is, the points $x
\in \bd D^2$ and $y \in \compb^1$ are identified on one hand if
and only if $y=f(\varphi_C(x))$, on the other hand, if and only if
$y=\varphi_{f(C)}(x)$, which yields that
$f(\varphi_C(x))=\varphi_{f(C)}(x)$ for all $x \in S^1$.

Regard $S^1$ as $[0,1]$ with its endpoints glued together to $x_0$
in such a way that $\varphi_{f(C)}(x_0)=\alpha(f(C)) \in \compb^0$.
Then for the paths corresponding to the attaching maps, we have $f
\circ p_C=p_{f(C)}$, that is, $f(\bound C)=\bound{f(C)}$. In
particular, $\alpha(f(C))=f(\alpha (C))$. Since $f$ respects the
labeling, this also yields $\ell (\bound C)= \ell (f(\bound C))=
\ell (\bound {f(C)})$.

To prove the uniqueness of $\bound C$, all we need to prove is that
$\varphi_C(x_0)$ is uniquely determined, as the label of the boundary walk of
$C$ and the root $\alpha(C)=\varphi_C(x_0)$ determine
$\bound C$ uniquely.
Take a neighborhood $N$ of $x_0$ in the disk $D^2$. Denote the
images of $N$ in $\comp$ and $\compb$  by $N_\comp$ and $N_\compb$
respectively after the identifications $ x \sim \varphi_C(x)$ and
$x \sim \varphi_{f(C)}(x)$ for $x \in \bd D^2$.
Naturally, $\varphi_C(x_0) \in N_\comp$ and $\varphi_{f(C)}(x_0) \in N_\compb$. Since
$f|_C$ is a homeomorphism, it takes $int(N_\comp)$ to
$int(N_\compb)$ homeomorphically, and therefore takes $N_\comp$ to
$N_\compb$. If $N$ is small enough, there is only one preimage of
$\varphi_{f(C)}(x_0)$ in $N_\compb$, and that is $\varphi_C(x_0)=\alpha(f(C))$.

\end{Proof}

\begin{Rem}

We point out that the second part of the theorem is non-trivial when
$\ell (\bound C)=x^n$ for some word $x$, in which case there may be
more than one vertex on $\bound C$ from which $\ell (\bound C)$ can
be read.

\end{Rem}

We have just seen that  for an immersion $f \colon {\comp}
\rightarrow {\compb}$ and for any $2$-cell $C \in {\comp}^2$, we
have $\ell (\bound C)=\ell (\bound{f(C)})$. In particular, when
$\compb = B_{X,P}$, then for any $2$-cells $C_1, C_2 \in \comp$ with
$\ell (C_1)=\ell (C_2)=\rho$, we have $\ell (\bound{C_1})=\ell
(\bound{C_2})$: this common label (called the ``boundary label" of
$\rho$) will often be denoted by $\bl (\rho)$. Thus $\bl (\rho) \in
(\X)^*$.

As in covering space theory, paths of a $2$-complex $\comp$ are our tools to classify
immersions over $\comp$. The point of the following construction is to generalize the
notion of graph-theoretic paths to $2$-complexes.

We associate an edge-labeled graph $\graph_\comp$ with the
2-complex $\comp$ as follows:
$$V(\graph_\comp)= \comp^{(0)}$$
$$E(\graph_\comp)= \comp^{(1)} \cup \{e_C: C \in \comp^{(2)}\},$$
 where $e_C$ denotes a loop based at $\init C$ and labeled by $\ell (C)$. Thus the edges
 in $\comp^{(1)}$ are labeled over $\X$ and the edges of the form $e_C$ (for $C$ a $2$-cell) are
 labeled over $P$.
Since an edge labeled by $\rho \in P$ is always a loop, we may identify $P$ with
 $P^{-1}$ and regard $\graph_\comp$ as an $X \cup P$-graph in the
 sense of section 2 of the paper.

\begin{Lem}
For any labeled complex $\comp$, the labeled graph $\graph_\comp$ is deterministic.
\end{Lem}

\begin{Proof}
Let $f \colon \comp \to B_{X,P}$ be the immersion inducing the
labeling on $\comp$. The subgraph corresponding to the $1$-skeleton
of $\comp$ is deterministic, as its labeling is induced by the
immersion $f|_{\comp^1}$ over $B_X$ (see \cite{MM1}). Therefore we
only need to check if different edges labeled by $\rho$ are based at
different vertices, that is, if different $2$-cells in $\comp$
labeled by $\rho$ have different roots. Denote the set of
$\rho$-labeled $2$-cells of $\comp$ by $\{C_\alpha: \alpha \in A\}$,
and the corresponding attaching maps $\varphi_\alpha \colon S^1 \to
\comp$ for $\alpha \in A$. Again, regard $S^1$ as the unit interval
with its endpoints identified with $x_0$, and let $N$ be a
neighborhood of $x_0$ in the disk $D^2$. Let $N_\alpha$ denote the
image of $N$ induced by the attaching map $\varphi_\alpha$. Since
$f$ maps all $\rho$-labeled $2$-cells to one cell,
$f(N_\alpha)=f(N_{\alpha'})$ for all $\alpha, \alpha' \in A$ and for any neighborhood $N$.
Since
$f$ is locally injective, this implies that the $f(N_\alpha)$
$(\alpha \in A)$ are pairwise disjoint, therefore the roots
$\varphi_\alpha(x_0)$ of the $2$-cells are all different.
\end{Proof}

The paths in the graph $\graph_\comp$ will play the role of paths in
$\comp$ in our paper --- we classify immersions by keeping track of how these paths lift. One can think of these paths as paths in
$\comp^1$ (in the graph-theoretic sense) extended with the
possibility of ``stepping" on a $2$-cell at its basepoint, thus
including it in the path. We remark that if the immersion is
actually a cover, then every $2$-cell in the base space lifts to a
$2$-cell based at every point in the fibre of a point in the base
space, so it is not necessary to keep track of the loops $e_C$ in
this case.

\begin{Lem}
\label{comp-graph-immersion} For two labeled $2$-complexes $\comp$
and $\compb$ there exists an immersion $\comp \to \compb$ (that
respects the labeling) if and only if there is an immersion
$\graph_\comp \to \graph_\compb$ (that respects the labeling).
\end{Lem}

\begin{Proof}
Let $f \colon \comp\to \compb$ be an immersion that respects the
labeling. Regarding $\comp^1$ as a subgraph of $\graph_\comp$, we
define $g \colon \graph_\comp \to \graph_\compb$ to be $f$ on
$\comp^1$, and for an edge $e_C$ corresponding to a $2$-cell $C$,
let $g(e_C)=e_{f(C)}$. It is easy to see that if $f$ is locally
injective at the vertices, so is $g$, hence an immersion. For the
converse, suppose $g \colon \graph_\comp \to \graph_\compb$ is an
immersion that respects the labeling. Define $f \colon \comp \to
\compb$ to be $g$ on $\comp^1$, and for a $2$-cell $C$ of $\comp$,
let $f(C)$ be the $2$-cell for which $g(e_C)=e_{f(C)}$ holds. Note
that if $g$ is an immersion, then so is $f|_{\comp^1}$. Suppose that
$f|_{\comp^1}$ is an immersion, but $f$ is not. Then there is a
vertex $v$ with two $2$-cells $C_1$ and $C_2$ with $v \in \bd C_1
\cap \bd C_2$ that $f$ identifies around $v$, that is, for any
neighborhood $N$ of $v$, $f(C_1 \cap N)=f(C_2 \cap N)$. Since $f$ is
locally injective on to the $1$-skeleton --- in particular, on
$\bound {C_1}$ and $\bound {C_2}$ ---, this can only happen if $C_1$
and $C_2$ have the same boundary walk, so $e_{C_1}$ and $e_{C_2}$
are based at the same vertex. But since $g(e_{C_1})=g(e_{C_2})$,
that contradicts our assumption. Hence $f$ is an immersion, and it
respects the labeling.

\end{Proof}

\section{Classification of immersions}

We are now ready to define an inverse monoid $\M{X,P}$ which will
play the role of the fundamental group. Let $\comp$ be a labeled
2-complex with the edges ($1$-cells) labeled over the set $\X$ and
the 2-cells labeled over the set $P$, consistent with an immersion
over some complex $B_{X,P}$. We  define the inverse monoid $\M{X,P}$
by

$$\M{X,P}=Inv \left\langle X \cup P\ |\ \rho^2=\rho,\ \rho \leq \bl (\rho) : \rho \in P \right\rangle,$$

or equivalently

$$\M{X,P}=Inv \left\langle X \cup P\ |\ \rho^2=\rho,\ \rho = \rho \, \bl (\rho) : \rho \in P \right\rangle$$

This monoid $\M{X,P}$ acts on the vertices ($0$-cells) of $\comp$
 as follows. For $x \in \X$, let $vx=w$ if there is an
edge labeled $x$ from $v$ to $w$, and $vx$ is undefined otherwise.
For $\rho \in P$, let $v\rho=v$ if there is a 2-cell labeled $\rho$
based at $v$, and $v\rho$ is undefined otherwise. This action
extends to an action of $\FIM(X)$ in a natural way. Since the action
of $\rho$ is always idempotent, and is always a restriction of the
action of $\bl (\rho)$, it also extends to an action of $\M{X, P}$.
Note that the action of $\M{X,P}$ on the vertices of $\comp$
corresponds to the usual partial action induced by edges in
$\graph_\comp$. We will denote the stabilizer of a vertex $v \in
{\comp}^0$ under this action by $\M{X,P}$ by $Stab({\comp},v)$.

\begin{Prop}
The inverse monoid $\M{X, P}$ and its action on $\comp^0$ do not
depend on the boundary walks and roots chosen for the $2$-cells.
\end{Prop}

\begin{Proof}
Suppose we chose different roots and boundary walks for the
$2$-cells of $\comp$, and let $\bl'(\rho)$ denote the new boundary
label corresponding to the $2$-cells labeled by $\rho$. The inverse
monoid corresponding to these boundary walks is $\M{X,P}'=\langle X,
P\ |\ \rho^2=\rho,\ \rho \leq  \bl' (\rho) \rangle$. The word
$\bl'({\rho})$ is a cyclic conjugate of $\bl(\rho)$ or
$(\bl(\rho))^{-1}$. Since $\rho = {\rho}^{-1}$, $\rho \leq \bl
(\rho)$ holds if and only if $\rho \leq (\bl (\rho))^{-1}$, so
reversing the boundary walk does not effect $\M{X,P}$. Hence  we may
assume that $\bl'({\rho})$ is a cyclic conjugate of $\bl(\rho)$.
Suppose $\bl(\rho)=p_\rho q_\rho$, $\bl'({\rho})=q_\rho p_\rho$.
Note that $p_\rho \rho p_{\rho}^{-1}$ is an idempotent of
$\M{X,P}'$, since $\rho$ is an idempotent of $\M{X,P}'$. Also, since
${\rho} \leq q_{\rho}p_{\rho}$ in $\M{X,P}'$, it follows that
$p_{\rho}{\rho}p_{\rho}^{-1} =
p_{\rho}{\rho}q_{\rho}p_{\rho}p_{\rho}^{-1} \leq
p_{\rho}{\rho}q_{\rho} \leq p_{\rho}q_{\rho}$ in $\M{X,P}'$. Hence
the map $x \mapsto x, \rho \mapsto p_\rho \rho p_\rho^{-1}$, where
$x \in X$, $\rho \in P$, extends to a well-defined morphism $\varphi
\colon \M{X,P} \to \M{X,P}'$. Also, for $\rho \in \M{X,P}'$,
$p_\rho^{-1}(p_\rho \rho p_\rho^{-1}) p_\rho=\rho$, so $\varphi$ is
surjective; and it is injective since it is injective on the
generators of $\M{X,P}$, so it is an isomorphism.

Moreover, denoting the maps from $\M{X, P}$ and $\M{X,P}'$ to
$\SIM(\comp^0)$ corresponding to their actions on the vertices by
$\psi$ and $\psi'$ respectively, the following diagram commutes:

\begin{picture}(160,95)(-70,-70)
\put(40,10){\circle*{3.5}} \put(144,10){\circle*{3.5}}
\put(92,-38){\circle*{3.5}}

\put(32,6){\makebox(0,0)[br]{$\M{X,P}$}}
\put(152,6){\makebox(0,0)[bl]{$\M{X,P}'$}}
\put(92,-56){\makebox(0,0)[bc]{$\SIM(\comp^0)$}}
\put(89,11){\makebox(0,0)[bl]{$\varphi$}}

\put(54,-20){\makebox(0,0)[br]{$\psi$}}
\put(130,-20){\makebox(0,0)[bl]{$\psi'$}}
\put(45,5){\vector(1,-1){42}}
\put(139,5){\vector(-1,-1){42}}
\put(46,9){\vector(1,0){93}}
\end{picture}

The commutativity of the diagram follows directly from the facts
that $\varphi$ is the identity on $X$, and that for $\rho \in \M{X,P}$, the action of
$\varphi(\rho)$ on the vertices is the same as that of $\rho$.
\end{Proof}

We now define an inverse category of paths on $\graph_\comp$. A
category $\cat$ is called \emph{inverse} if for every morphism $p$
in $\cat$ there is a unique inverse morphism $p^{-1}$ such that
$p=pp^{-1}p$ and $p^{-1}=p^{-1}pp^{-1}$. The loop monoids $L(\cat,
v)$ of an inverse category, that is, the set of all morphisms from
$v$ to $v$, where $v$ is an arbitrary vertex, form an inverse
monoid. The free inverse category $\operatorname{FIC}(\graph)$ on a
graph $\graph$ is the free category on $\graph$ factored by the
congruence induced by relations of the form $p=pp^{-1}p$,
$p^{-1}=p^{-1}pp^{-1}$, and $pp^{-1}qq^{-1}=qq^{-1}pp^{-1}$ for all
paths $p, q$ in $\graph$ with $\init p=\init q$.

Now let $\sim$ be the congruence on the free category on
$\graph_\comp$ generated by the relations defining
$\operatorname{FIC}(\graph_\comp)$ and the ones of the form $p^2=p$
and $p = pq$, where $p, q$ are coterminal paths with $\ell (p) \in
P$ and $\ell (q)=\bl (\ell(p))$. The inverse category $\pcat \comp$
corresponding to the $2$-complex $\comp$ is obtained by factoring
the free category on $\graph_\comp$ by $\sim$. The loop monoids $L(\pcat
\comp, v)$ consist of $\sim$-classes of $(v,v)$-paths, these monoids
play the role of the fundamental group, and $\pcat \comp$ plays the
role of the fundamental groupoid in the classification of
immersions. We will denote $L(\pcat \comp, v)$ by $L(\comp, v)$ for
brevity.

\begin{Prop}
\label{fundgroup}
For any vertex $v$ in a connected $2$-complex
$\comp$, the greatest group homomorphic image of $L(\comp, v)$ is
the fundamental group of $\comp$.
\end{Prop}

\begin{Proof}
The proof follows from the fact that the fundamental groupoid of
$\comp$ is $\pcat \comp$ factored by the congruence generated by
relations of the form $xx^{-1}=\operatorname{id}_{\alpha(x)} $ for
any morphism $x$ (which implies $\bound C=\operatorname
{id}_{\alpha(C)}$ for any $2$-cell $C$). Hence $L(\comp,
v)/\sigma=\pi_1(\comp)$.
\end{Proof}

Note that the relations of $\sim$ are closely related to the ones
defining $\M{X,P}$, that is, two coterminal paths $p, q$ are in the
same $\sim$-class if and only if $\ell (p)=\ell (q)$ in $\M{X,P}$.
This enables us to identify morphisms from some vertex $v$ with
their (common) label in $\M{X,P}$. Using this identification, we
have $L(\comp, v)=\Stab(\comp, v)$ for any vertex $v$. Note that if $\comp =
B_{X,P}$, then $\pcat \comp = \M{X,P}$. The following
proposition is a direct consequence of our previous observation and
Proposition \ref{stabclosed}.

\begin{Prop}
\label{loopclosed} Each loop monoid of $\pcat \comp$ is a closed
inverse submonoid of $\M{X,P}$.
\end{Prop}

Given a closed inverse submonoid $H$ of $\M{X,P}$, we construct a complex with $H$ as a
loop monoid using the $\omega$-coset graph $\graph_H$ of $H$.
First note that the action of
$\M{X, P}$ by right multiplication on the right $\omega$-cosets of
$H$ is by definition the same as the action on the vertices of
$\graph_H$ induced by the edges. Suppose there is a closed path
based at $H$ labeled by $x\rho y$, where $\rho \in P,\ x,y \in (X
\cup X^{-1} \cup P)^\ast$. Then $x\rho y \in H$, and since $H$ is
closed and $x \rho y \leq xy$ in $\M{X, P}$, we also have $xy \in
H$, hence $xy$ also labels a closed path based at $H$. This implies
that $\rho$ always labels a loop in the coset graph. Similarly, $x
\rho y \leq x (\bl (\rho)) y$, so $x (\bl (\rho)) y$ labels a closed
path based at $H$. Therefore whenever there is a loop in the coset
graph labeled $\rho$ based at $v$, there is a closed path labeled
$\bl (\rho)$ based at $v$.

The labeled \textsl{coset complex} $\comp_H$ of $H$ is defined the
following way:
$$\comp_H^{(0)}=V(\graph_H),$$
$$\comp_H^{(1)}=\{e \in E(\graph_H): \ell (e) \in X \cup X^{-1}\},$$
$$\comp_H^{(2)}=\{C_e \in E(\graph_H): \ell (e) \in P\},$$
where the boundary walk of a 2-cell $C_e$ is the closed path rooted
at $\init e$ and labeled  by $\bl ({\rho})$ where $\rho = \ell (e)$.
In short, we take the graph $\graph_H$, and substitute edges labeled
by $P$ with 2-cells in the natural way. Note that the labeling of
$C_H$ corresponds to the immersion over the $2$-complex $B_{X,P}$,
in which the attaching map of a $2$-cell labeled by $\rho$ is given
by $\bl (\rho)$.

The following proposition gives the relationships between the
complexes associated with the coset graphs and graphs associated with complexes.

\begin{Prop}
\label{stabcoset} Let $\comp$ be a labeled $2$-complex. If $H$ is a
closed inverse submonoid of $\M{X, P}$ for which $\graph_H \cong
\graph_\comp$, then $\comp_H \cong \comp$. There is an isomorphism
$\varphi \colon \graph_H \to \graph_\comp$ if and only if
$H=\Stab(\comp,v)$ for some vertex $v$, and in that case, $\varphi$ is uniquely determined and $\varphi(H)=v$.

\end{Prop}

\begin{Proof}
The first statement follows directly from the definitions of
$\comp_H$ and $\graph_\comp$. For the second statement, suppose $H =
\Stab(\comp, v)$ for some $v \in \comp^0$. First we observe that the
set of words labeling closed paths from $\Stab(\comp, v)$ to
$\Stab(\comp, v)$ in $\graph_{\Stab(\comp, v)}$ is the same as the
set of words labeling closed paths from $v$ to $v$ in
$\graph_\comp$. Indeed, $p$ is a closed $(v,v)$-path in
$\graph_\comp$ if and only if $\ell (p) \in \Stab(\comp, v)$, which
is if and only if $p$ is a closed path from $\Stab(\comp, v)$ to
$\Stab(\comp, v)$ in $\graph_{\Stab(\comp, v)}$. We now define an
isomorphism $\varphi \colon \graph_{\Stab(\comp, v)} \to
\graph_\comp$ by $\Stab(\comp, v) \mapsto v$, and all $(\Stab(\comp,
v),\Stab(\comp, v))$-paths map to the (unique) $(v,v)$-path with the
same label. It is routine to verify that this is a graph
isomorphism.

Now for the converse, suppose $H \neq \Stab(\comp, v)$ for any
vertex $v$. Then the set of labels of closed $(H,H)$-paths in
$\graph_H$ and the ones of closed $(v,v)$ paths in $\graph_\comp$
are different, for all $v \in V(\graph_\comp)$, hence the two graphs
cannot be isomorphic.
\end{Proof}

Let $H,K$ be two closed inverse submonoids of $M_{X,P}$. Define $H$
to be \emph{conjugate} to $K$, denoted by $H \approx K$, if there
exists $m \in M_{X,P}$ such that $m^{-1}Hm \subseteq K$ and
$mKm^{-1} \subseteq H$. It is easy to see that $\approx$ is an
equivalence relation (called ``conjugation") on the set of closed
inverse submonoids of $M_{X,P}$. The equivalence classes of
$\approx$ are called \emph{conjugacy classes}. We remark that
conjugate closed inverse submonoids of $\M{X,P}$ are not necessarily isomorphic (see \cite{MM1}).

We call the two (labeled) immersions $f_1 \colon \comp_1 \to \compb$ and
$f_2 \colon \comp_2 \to \compb$ \emph{equivalent} if there is a labeled
isomorphism $\varphi \colon \comp_1 \to \comp_2$ which makes the following diagram commute:

\begin{picture}(160,95)(-70,-70)
\put(40,10){\circle*{3.5}} \put(144,10){\circle*{3.5}}
\put(92,-38){\circle*{3.5}}

\put(32,6){\makebox(0,0)[br]{$\comp_1$}}
\put(152,6){\makebox(0,0)[bl]{$\comp_2$}}
\put(92,-56){\makebox(0,0)[bc]{$\compb$}}
\put(89,11){\makebox(0,0)[bl]{$\varphi$}}

\put(54,-20){\makebox(0,0)[br]{$f_1$}}
\put(130,-20){\makebox(0,0)[bl]{$f_2$}}
\put(45,5){\vector(1,-1){42}}
\put(139,5){\vector(-1,-1){42}}
\put(46,9){\vector(1,0){93}}
\end{picture}

The following two theorems state the main result of the paper. They are generalizations of
Theorem 4.4 and
4.5 in \cite{MM1}, and most of the proofs are  analogous to those. When the
$2$-complexes contain no $2$-cells (that is, they are graphs), these theorems
reduce to Theorem 4.4 and
4.5 in \cite{MM1}.

\begin{Thm}
\label{premain} Let $\comp$ be a 2-complex, with edges labeled over
the set $\X$, 2-cells labeled over the set $P$, consistent with an
immersion over some complex $B_{X, P}$. Then each loop monoid is a
closed inverse submonoid of $\M{X,P}$, and the set of all loop
monoids $L(\comp,v)$ for $v \in \comp^0$ forms a conjugacy class of
the set of closed inverse submonoids of $\M{X, P}$. Conversely, if
$H$ is a closed inverse submonoid of $\M{X, P}$, then there is a
2-complex $\comp$ and an immersion $f \colon \comp \to B_{X, P}$
such that $H$ is a loop monoid of $\pcat\comp$, furthermore, $\comp$
is unique (up to isomorphism), and $f$ is unique (up to
equivalence).
\end{Thm}

\begin{Proof}
We saw in Proposition \ref{loopclosed} that loop monoids are closed.
Take two loop monoids $L(\comp, v_1)$ and $L(\comp, v_2)$, and let
$m \in (X \cup X^{-1} \cup P)^*$ label a $(v_1,v_2)$-path in
$\comp$. If $n \in L(\comp, v_2)$, then $n$ labels a $(v_2,
v_2)$-path, and $mnm^{-1}$ labels a $(v_1,v_1)$-path, so $m L(\comp,
v_2)m^{-1} \subseteq L(\comp, v_1)$. Since $m^{-1}$ labels a
$(v_2,v_1)$-path, we get $m^{-1} L(\comp, v_1)m \subseteq L(\comp,
v_2)$ similarly. Now suppose $H \approx L(\comp, v_1)$. Then there
exists some $m \in \M{X,P}$ such that $m^{-1}L(\comp, v_1)m=H$ and
$mHm^{-1} =L(\comp, v_1)$, in particular, $mm^{-1} \in L(\comp,
v_1)$. Therefore, regarding $m$ as an element of $(X \cup X^{-1}
\cup P)^\ast$, it labels a path from $v_1$ to some vertex $v_2$. If
$h \in H$ (and again regard $h$ as an element of $(X \cup X^{-1}
\cup P)^\ast$), then $mhm^{-1}$ labels a $(v_1,v_1)$-path, hence $h$
labels a path form $v_2$ to $v_2$. Therefore $H \subseteq L(\comp,
v_2)$. On the other hand, if $n \in L(\comp, v_2)$, then $mnm^{-1}
\in L(\comp, v_1)$, and $m^{-1}mnm^{-1}m \subseteq H$. Since $H$ is
closed and $m^{-1}mnm^{-1}m \leq n$, this yields $n \in H$,
therefore $H=L(\comp, v_2)$. This proves that the set of all loop
monoids $L(\comp,v)$ for $v \in \comp^0$ form a conjugacy class of
the set of closed inverse submonoids of $\M{X, P}$.

Now suppose that $H$ is a closed inverse submonoid of $\M{X, P}$,
and build the coset complex $\comp_H$. There is a natural immersion
$f\colon \comp_H \to B_{X, P}$, namely the one sending all edges and
2-cells to the ones corresponding to their labels.

It follows from Proposition \ref{stabcoset} that the graph
$\graph_\comp$ is unique, and it uniquely determines $\comp$. The
uniqueness of $f$ follows from the fact that $f$ respects the
labeling.

\end{Proof}

\begin{Thm}
\label{main}
Let $f \colon \comp_2 \to \comp_1$ be an immersion over $\comp_1$,
where $\comp_1$ and $\comp_2$ are 2-complexes with edges labeled
over the set $X \cup X^{-1}$, 2-cells labeled over the set $P$
consistent with an immersion over some complex $B_{X, P}$, and $f$
respects the labeling. If $v_i \in \comp_i^0$, $i=1,2$, such that
$f(v_2)=v_1$, then $f$ induces an embedding of $L(\comp_2, v_2)$
into $L(\comp_1, v_1)$. Conversely, let $\comp_1$ be a labeled
2-complex  and let $H$ be a closed inverse submonoid of $\M{X,P}$
such that $H \subseteq L(\comp_1, v_1)$ for some $v_1 \in
\comp_1^0$. Then there exists a 2-complex $\comp_2$ and an immersion
$f \colon \comp_2 \to \comp_1$ and a vertex $v_2 \in \comp_2^0$ such
that $f(v_2)=v_1$ and $L(\comp_2, v_2)=H$. Furthermore, $\comp_2$
is unique (up to isomorphism), and $f$ is unique (up to
equivalence). If $H, K$ are two closed inverse submonoids of
$\M{X,P}$ with $H, K \subseteq L(\comp_1, v_1)$, then the
corresponding immersions are equivalent if and only if $H \approx K$
in $\M{X,P}$.
\end{Thm}

\begin{Proof}
Suppose first that $f(v_{2}) = v_1$. The assertion that $f$ induces an embedding from $L(\comp_2,
v_2)$ to $L(\comp_1, v_1)$ follows easily from the fact that if
$p$ is a closed path in $\comp_2$ based at $v_2$, then $f(p)$ is a
closed path in $\comp_1$ based at $f(v_2)=v_1$ and $\ell (p) = \ell
(f(p))$. For the converse, suppose $H$ is a closed inverse submonoid
of $\M{X,P}$ such that $H \subseteq L(\comp_1, v_1)$, and construct
the coset complex $\comp_H$ and the coset graph $\graph_H$, and let
$\graph_1$ denote $\graph_{\comp_1}$. Put $\comp_2=\comp_H$, and
$v_2=H$. We saw in Proposition \ref{stabcoset} that $H = L(\comp_H,
H)$. We construct an immersion $g \colon \graph_H \to \graph_1$ that
respects the labeling. Let $f(H)=v_1$, and note that if
$(Hm)^{\omega}$ is a right $\omega$-coset, then $mm^{-1} \in H
\subseteq L(\comp_1, v_1)$, so $m$ labels a path starting at $v_1$
in $\graph_1$. Now we define $g$ to take all paths starting at $H$
to the (unique) path with the same label, starting at $v_1$. Then
$g$ is locally injective at the vertices, hence it is an immersion,
and it respects the labeling by definition. By Lemma
\ref{comp-graph-immersion}, $g$ yields an immersion $f \colon
\comp_H \to \comp_1$ that commutes with the labeling.

The uniqueness
of $f$ and $\comp_2$ again follow from the uniqueness of
$\graph_{\comp_2}$ by Proposition \ref{stabcoset}, and from the fact
that $f$ respects the labeling. For the last statement, recall that
according to Lemma \ref{attaching}, the immersion $f$ and the
complex $\comp_2$ determine the boundary walks and therefore the
graph $\graph_{\comp_2}$ uniquely, that is, $\graph_{\comp_2}$ and
the pair $(f, \comp_2)$ are in one-one correspondence. The fact
$\graph_H \cong \graph_{H'}$ if and only if $H$ and $H'$ are
conjugate completes our proof.
\end{Proof}

We close this section with some observations about the inverse
monoids $\M{X, P}$ and their closed inverse submonoids. In
particular, we give an algorithm to construct $\comp_H$ for a
finitely generated closed inverse submonoid $H$ of $\M{X,P}$ if $X$
and $P$ are finite.

\begin{Thm}
\label{mainalg}
(a) If $X$ and $P$ are finite sets, then the
Sch\"utzenberger graphs of $\M{X, P}$ are finite (and effectively
constructible) and so the word problem for $\M{X, P}$ is decidable.

(b) If $X$ and $P$ are finite sets and $H$ is a finitely generated
closed inverse submonoid of $\M{X, P}$, then the  associated
$2$-complex $\comp$ is finite and effectively constructible.

\end{Thm}

\begin{Proof} (a) If $w$ is a word in $(\X)^*$ then no defining
relation for $\M{X,P}$ applies, so the corresponding
Sch\"utzenberger graph $S{\graph}(w)$ is the Munn tree of $w$ (see
\cite{munn, law}), so it is finite and effectively constructible. On
the other hand, if $w$ is a word in $(\X \cup P)^*$ that does
contain some letter $\rho \in P$, then any application of the
relation ${\rho}^{2} = \rho$ turns the edge labeled by $\rho$ into a
loop. Any application of the relation $\rho \leq \bl (\rho)$ (i.e.
$\rho = \rho \, \bl (\rho)$) just introduces a new path labeled by
$\bl (\rho)$ to the approximate automaton. Once the relations $\rho
= {\rho}^2$ and $\rho \leq \bl (\rho)$ have been applied, this
occurrence of $\rho$ is not involved in any further application of
relations involved in iteratively constructing $S{\graph}(w)$. As
the automaton we started out with was finite, this iterative process
(as outlined in Section 2 above - Theorem 4.12 of \cite{Ste1}) must
terminate in a finite number of steps and the Sch\"utzenberger
automaton $S{\mathcal A}(w)$ is finite and effectively
constructible.

(b) The proof of part (b) of the theorem is similar. If we start
with the flower automaton ${\mathcal F}(Y)$ of a finite subset $Y
\subset (\X \cup P)^*$ and iteratively apply edge foldings and
expansions corresponding to the defining relations of $\M{X,P}$,
this process terminates in a finite number of steps, providing an
effective construction of the
 $\omega$-coset automaton of the corresponding closed inverse
submonoid $\langle Y\rangle ^{\omega} $ of $\M{X,P}$ by Theorem
\ref{cosetalg}. The result then follows from Theorem 4.8 (and the
fact that the associated complex $\comp$ is the coset complex of
$\langle Y\rangle ^{\omega} $).

\end{Proof}

\section{Examples and special cases}

Recall that a \emph{covering space} of a space $X$ is a space
$\tilde X$ together with a map $f \colon \tilde X \to X$ called a
covering map, satisfying the following condition: there exists an
open cover ${U_\alpha}$ of $X$ such that for each $\alpha$,
$f^{-1}(U_\alpha)$ is a disjoint union of open sets in $\tilde X$,
each of which is mapped homeomorphically onto $U_\alpha$ by $f$. It
is easy to see  that a cellular  map $f \colon \comp \to \compb$
between CW-complexes is a covering map if and only if each $0$-cell
$v \in C^0$ has a neighborhood $U_v$ that is homeomorphic to a
neighborhood $U_{f(v)}$ of $f(v)$. This happens if and only if $f$
is an immersion for which the neighborhoods of $0$-cells ``lift
completely", that is, whenever $v$ is on the boundary of a cell $C$
in $\compb$, then each $0$-cell in $f^{-1}(v)$ is on the boundary of
a cell in $f^{-1}(C)$.

The following theorem characterizes those immersions between
$2$-complexes which are also covering maps, in the sense of the
previous theorem.

\begin{Thm}
\label{cover} Let $\comp, \compb$ be $2$-complexes labeled by an
immersion over some complex $B_{X,P}$, let $f \colon \comp \to
\compb$ be an immersion that respects the labeling, and let $v \in
\comp^0$ be an arbitrary $0$-cell. Then $f$ is a covering map if and
only if
$L(\comp, v)$ is a full closed inverse submonoid of $L(\compb, f(v))$, that is, it
contains all idempotents of $L(\compb, f(v))$.
\end{Thm}

\begin{Proof}
First, suppose that $f$ is a covering,
and suppose there is an idempotent $e \in
L(\compb, f(v))$. Regarding $e$ as an element of $(X \cup X^{-1}
\cup P)^\ast$, the closed path in $\compb$ labeled by $e$, starting at
$f(v)$ lifts to a path labeled by $e$, starting at $v$ in $\comp$,
because $f$ is a covering. Since $e$ is idempotent, the action of
any path labeled by $e$ on $\comp^{0}$ is the restriction of the
identity, therefore a path labeled by $e$ must always be closed.
This yields $e \in L(\comp, v)$.

For the converse, suppose $L(\comp, v)$ is a full closed inverse submonoid of $L(\compb, f(v))$.
Suppose there is an edge starting at $f(v)$, labeled by $s$ in $\compb$. Then
$ss^{-1} \in L(\compb, f(v))$, and since $ss^{-1}$ is idempotent, that implies
$ss^{-1} \in L(\comp, v)$. Which yields that there is an edge labeled by $s$,
starting form $v$ in $L(\comp, v)$, that is, the
neighborhood of $f(v)$ lifts completely. By induction on distance from $v$, we
obtain that all $0$-cells are in the image of $f$, therefore their neighborhoods lift completely.

\end{Proof}

It is easy to see that $L(\comp, v)$ is a full closed inverse submonoid of $L(\compb, f(v))$
if and only if
whenever $m \in L(\comp, v)$ and $n \in L(\compb, f(v))$
such that $m \geq n$ holds in $L(\compb, f(v))$, then $n \in
L(\comp, v)$.
Therefore combining the result above with Theorem \ref{main}, we obtain that
an immersion $f \colon \comp \to \compb$ is a covering if and only
if whenever an element $m \in
L(\compb,f(v))$ is comparable with $n \in L(\comp, v)$ in the
natural partial order, we have $n \in L(\compb,f(v))$.

We briefly compare our results with the theorem classifying covers
via subgroups of the fundamental group when applied to
$2$-complexes. Recall (Proposition \ref{fundgroup}) that the
fundamental group $\pi_1(\comp)$ of a (connected) $2$-complex
$\comp$ is the greatest group homomorphic image of any loop monoid
of $\comp$, denoted by $L(\comp,v)/\sigma$. The greatest group
homomorphic image of $M_{X,P}$, denoted by $G_{X,P}$, is the group
with the same presentation as $\M{X,P}$. Since in groups,
$\rho=\rho^2$ implies $\rho=1$, and $\rho \leq \bl( \rho)$ implies
$\rho = \bl( \rho)=1$, that is just
$$G_{X,P}=Gp \langle X \ |\ \bl( \rho)=1  \rangle.$$
This is the fundamental group of the corresponding complex
$B_{X,P}$, and the fundamental group of a complex immersing into
$B_{X,P}$ is a subgroup of $G_{X,P}$.



Naturally, the fundamental groups of $2$-complexes immersing into a
$2$-complex $\comp$ are always subgroups of $\pi_1(\comp)$, but
distinct immersing $2$-complexes may give rise to the same subgroup
of $\pi_1(\comp)$
--- for example, any immersing tree has the trivial group as its
fundamental group. When restricting to covers, however, it is
well-known that the fundamental groups of the covering spaces are in
one-to-one correspondence with the conjugacy classes of subgroups of
the fundamental group of the base space. Therefore the loop monoids
of different covering spaces all have different greatest group
homomorphic images. Suppose $f \colon \comp \to \compb$ is a
covering that respects the labeling, and let $v \in \comp^0$. Let
$\sigma^\natural \colon L(\compb,f(v)) \to \pi_1(\compb)$ be the
natural homomorphism corresponding to the congruence $\sigma$.
Recall (\cite{law}) that $\sigma$ is generated by pairs $(m,n)$ such
that $m \leq n$. Therefore by Theorem \ref{cover} (and the
observation that followed), it is clear that $L(\comp, v)$ is the
union of some $\sigma$-classes of $L(\compb,f(v))$, namely it is the
full inverse image of $\pi_1(\comp)$ under $\sigma^\natural$.

In \cite{will}, Williamson uses similar methods to classify
immersions over a slightly restricted class of complexes with one
$0$-cell. The notion of immersion $f : \comp \rightarrow \mathcal D$
in \cite{will} has the additional property that every $0$-cell in
the fiber $f^{-1}(v_{0})$ of a $0$-cell $v_0$ on the boundary of a
$2$-cell of $\mathcal D$ is required to be part of the boundary of
some $2$-cell of $\comp$.

\begin{Ex}

Let $X=\{a\}$, $P=\{\rho\}$, and $\comp$ be the labeled $2$-complex with one loop labeled by
$a$ and one $2$-cell attached to the path $a^2$. (This complex is homeomorphic to the
projective plane.)
Then its loop monoid is $\M{X,P}=Inv\langle a, \rho\ |\ \rho^2=\rho,\ \rho=\rho a^2 \rangle$.
Here is a list of all $2$-complexes immersing into $\comp$, and a representative from the
corresponding conjugacy class of closed inverse submonoids of $\M{X,P}$. (The basepoint of
the representative is denoted by a larger dot when necessary.) The complex that immerses into
the projective plane uniquely determines the immersion (up to equivalence).
\vspace{10 pt}

\begin{longtable}{>{\centering}m{0.6\textwidth}<{\centering} | m{0.4\textwidth}}

\hline
\includegraphics[trim=0 0 0 -5, width=0.166\textwidth]{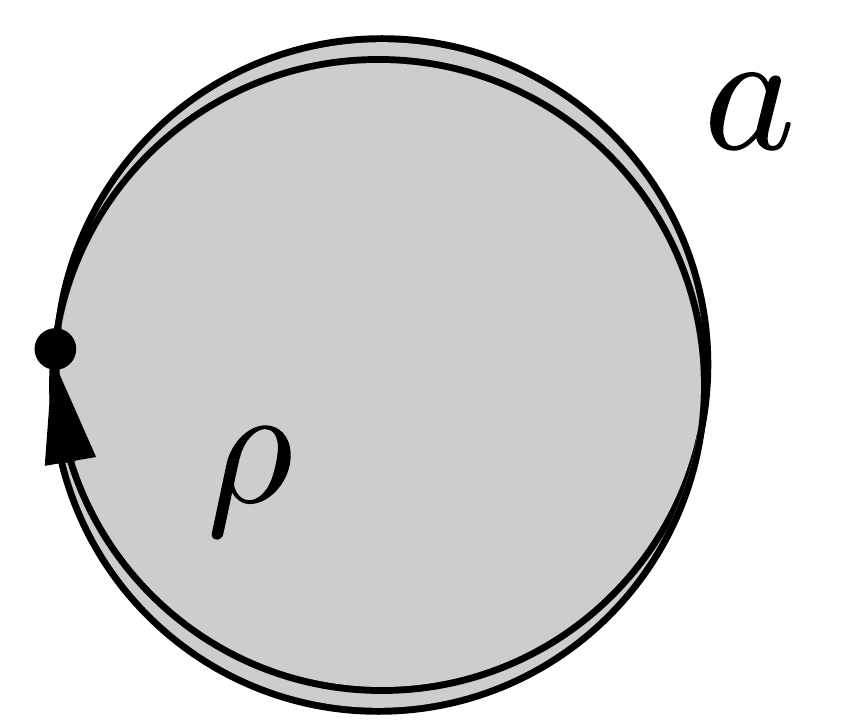}& $\langle a, \rho \rangle^\omega$, the projective plane\\
\hline
\includegraphics[trim=0 0 0 -5, width=0.20\textwidth]{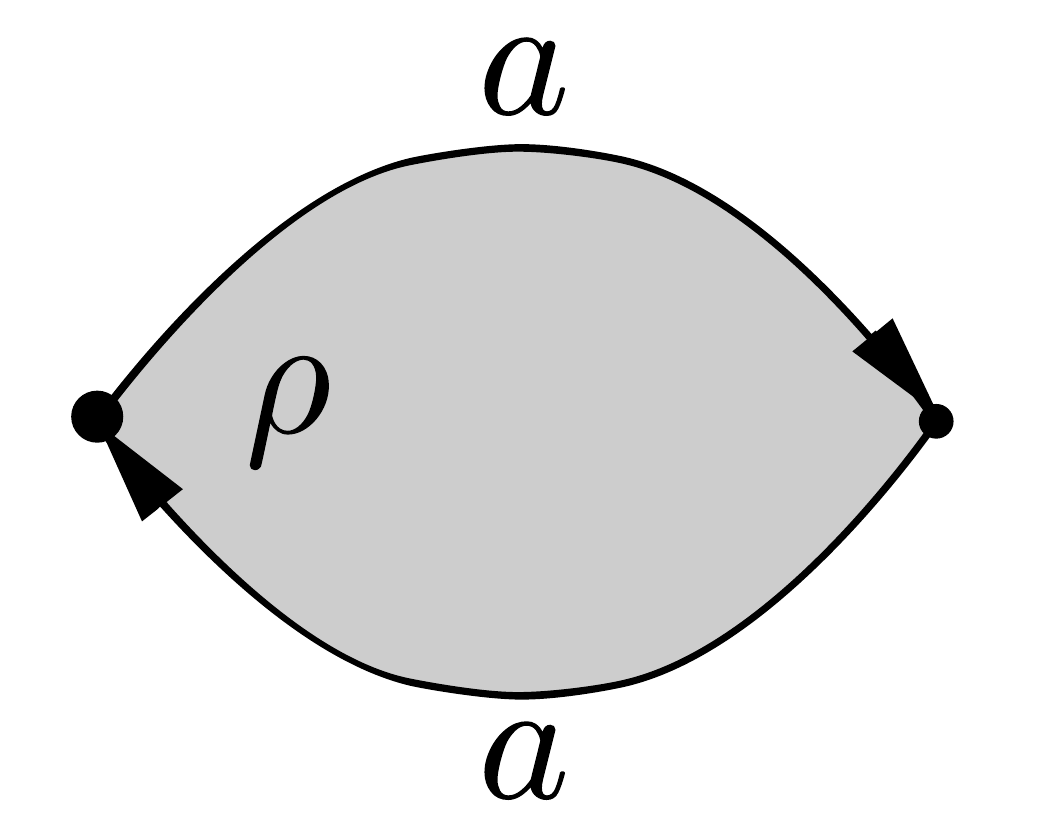}& $\langle \rho \rangle^\omega$\\
\hline
\includegraphics[trim=0 0 0 -5, width=0.20\textwidth]{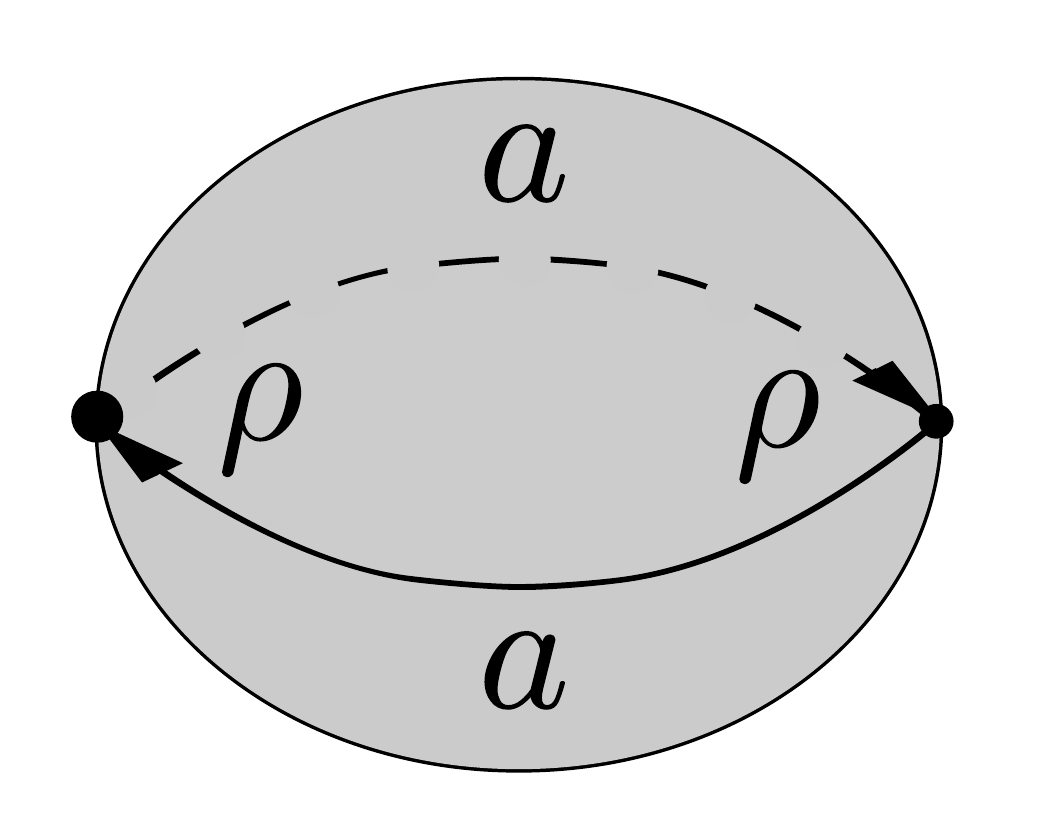}& $\langle \rho,\ a\rho a \rangle^\omega$, the universal cover\\
\hline
\includegraphics[trim=0 0 0 -5, width=0.166\textwidth]{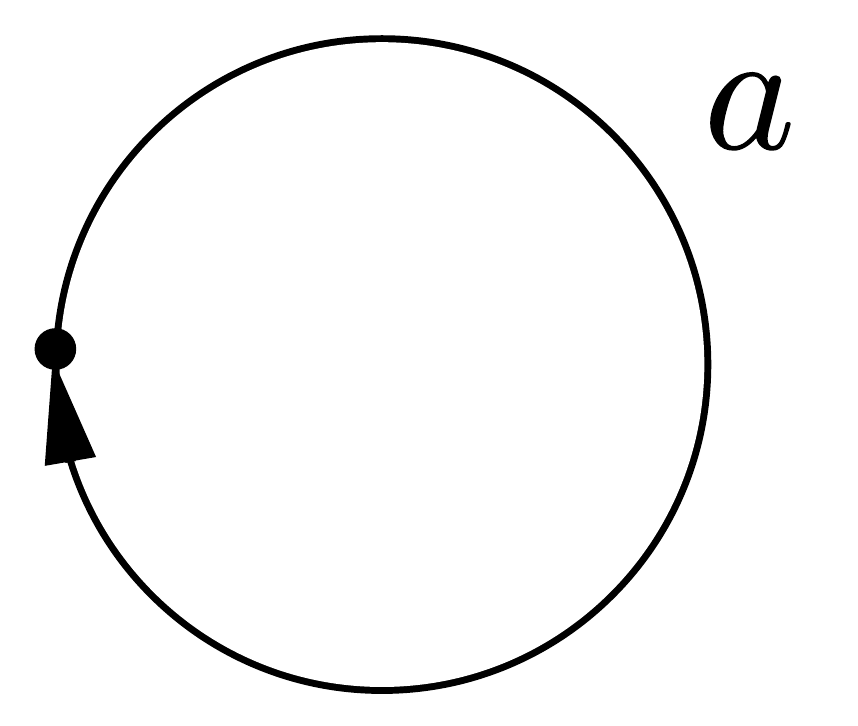}& $\langle a \rangle^\omega$\\
\hline
\includegraphics[trim=0 0 0 -5, width=0.30\textwidth]{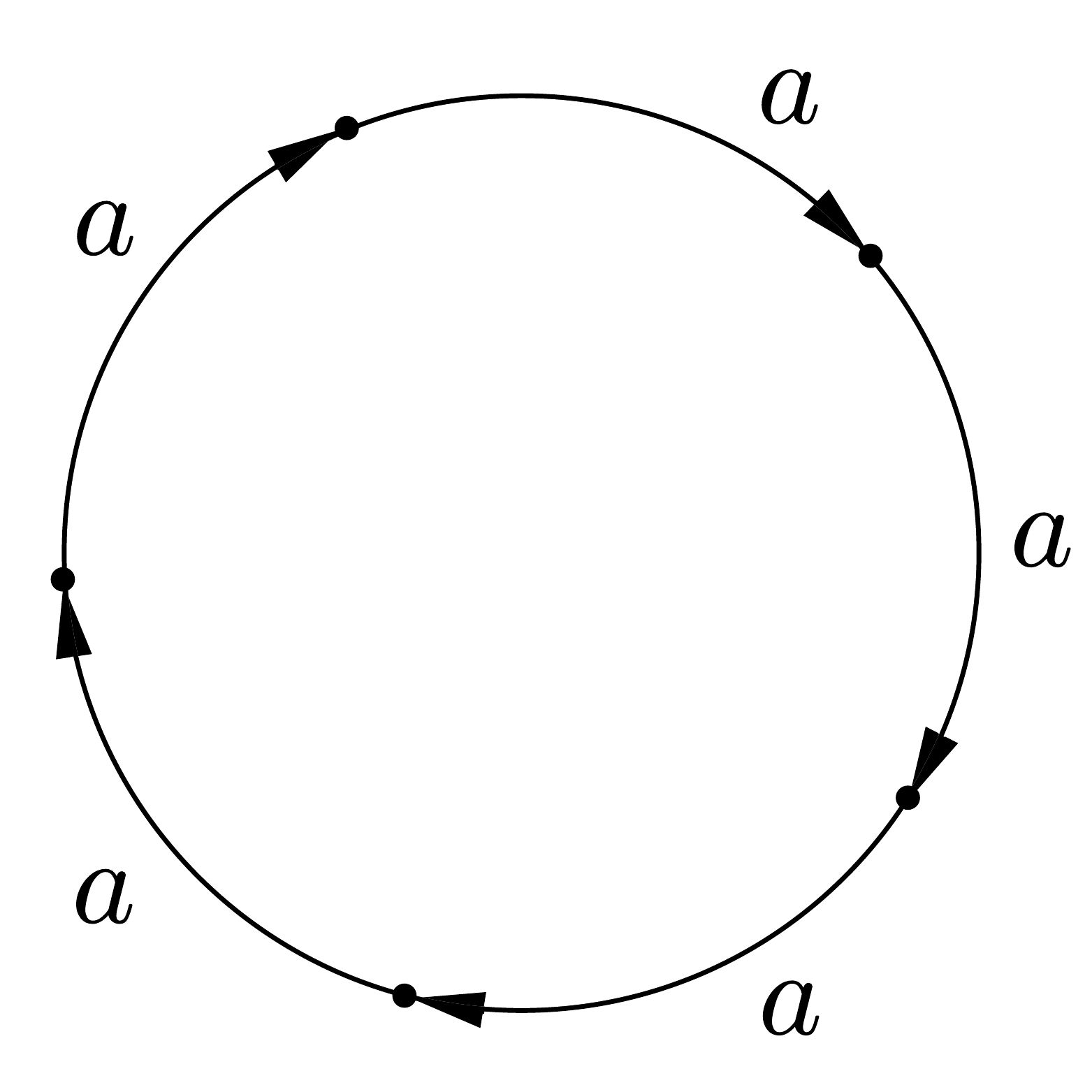}& $\langle a^n \rangle^\omega$, $n \in \mathbb N$, $(n=5)$\\
\hline
\includegraphics[width=0.02\textwidth]{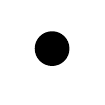}& $\langle 1 \rangle^\omega$\\
\hline
\includegraphics[trim=0 0 0 -5, width=0.40\textwidth]{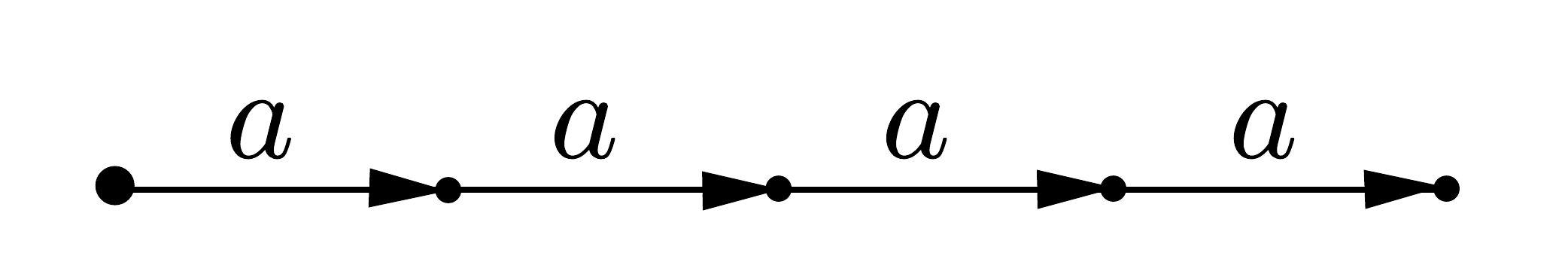}& $\langle a^na^{-n} \rangle^\omega$, $n \in \mathbb N$, $(n=4)$\\
\hline
\includegraphics[trim=0 0 0 -5, width=0.40\textwidth]{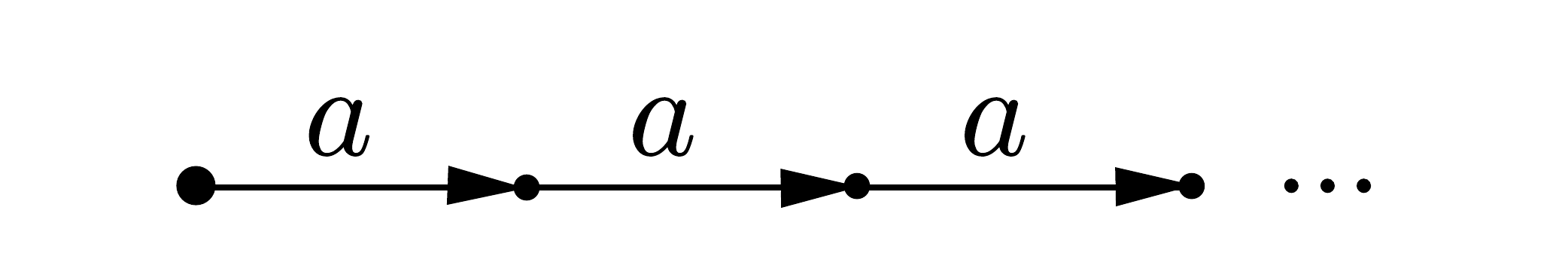}& $\langle a^na^{-n} : n \in \mathbb N \rangle^\omega$\\
\hline
\includegraphics[trim=0 0 0 -5, width=0.40\textwidth]{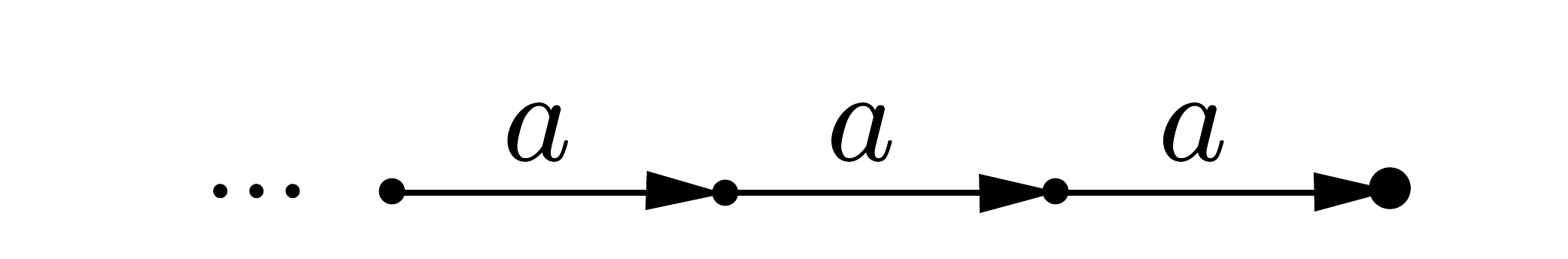}& $\langle a^{-n}a^{n} : n \in \mathbb N \rangle^\omega$\\
\hline
\includegraphics[trim=0 0 0 -5, width=0.40\textwidth]{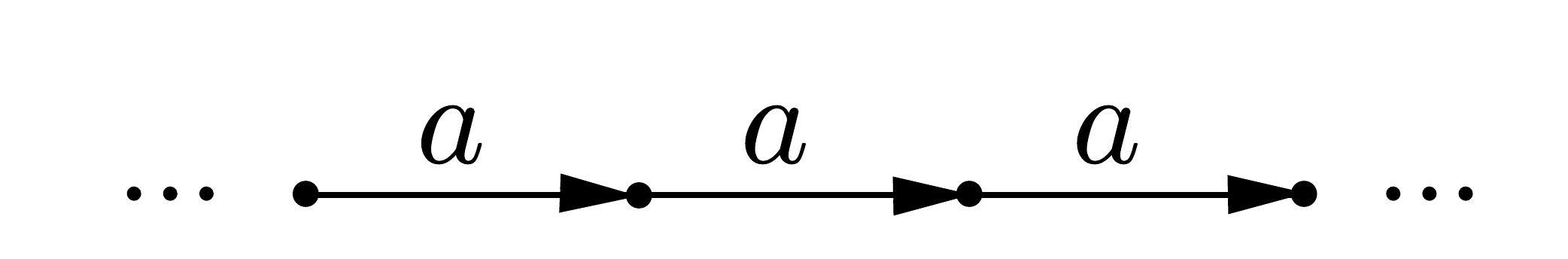}& $\langle a^{n}a^{-n} : n \in \mathbb Z \rangle^\omega$\\
\hline
\end{longtable}
\end{Ex}

\begin{Ex}

Let $X=\{a,b\}$, $P=\{\rho\}$, and let $\compb$ be the labeled
$2$-complex with two loops labeled by $a$ and $b$, and one $2$-cell
attached to the path $b$. Thus $\compb$ is a wedge sum of a circle
and a closed disk. Then its loop monoid is $\M{X,P}=Inv\langle
a,b,\rho \ |\ \rho^2=\rho,\ \rho=\rho b \rangle$. Here are some
examples of $2$-complexes immersing into $\compb$, and a
representative from the corresponding conjugacy class of closed
inverse submonoids of $\M{X,P}$.

\begin{longtable}{>{\centering}m{0.6\textwidth}<{\centering} | m{0.4\textwidth}}
\hline
\includegraphics[trim=0 0 0 -5, width=0.333\textwidth]{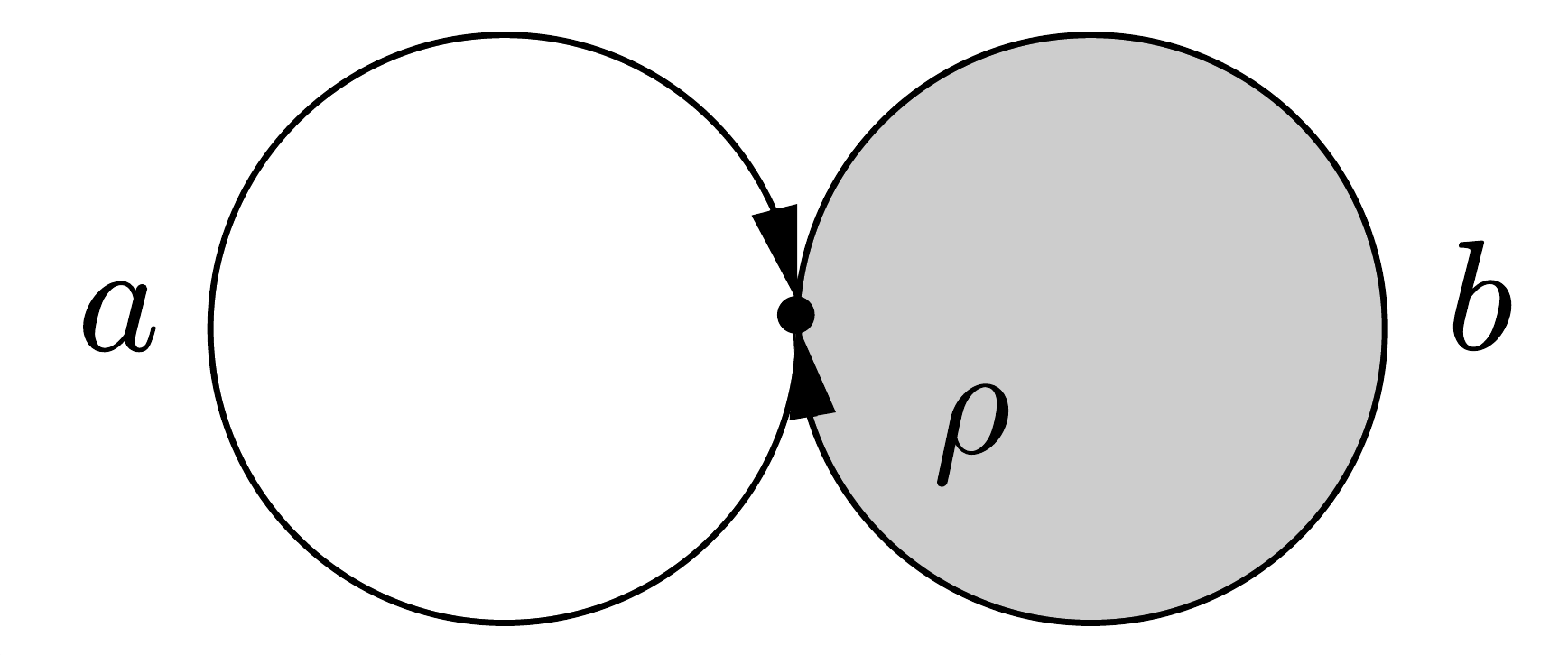}& $\langle a,b, \rho \rangle^\omega$\\
\hline
\includegraphics[trim=0 0 0 -5, width=0.333\textwidth]{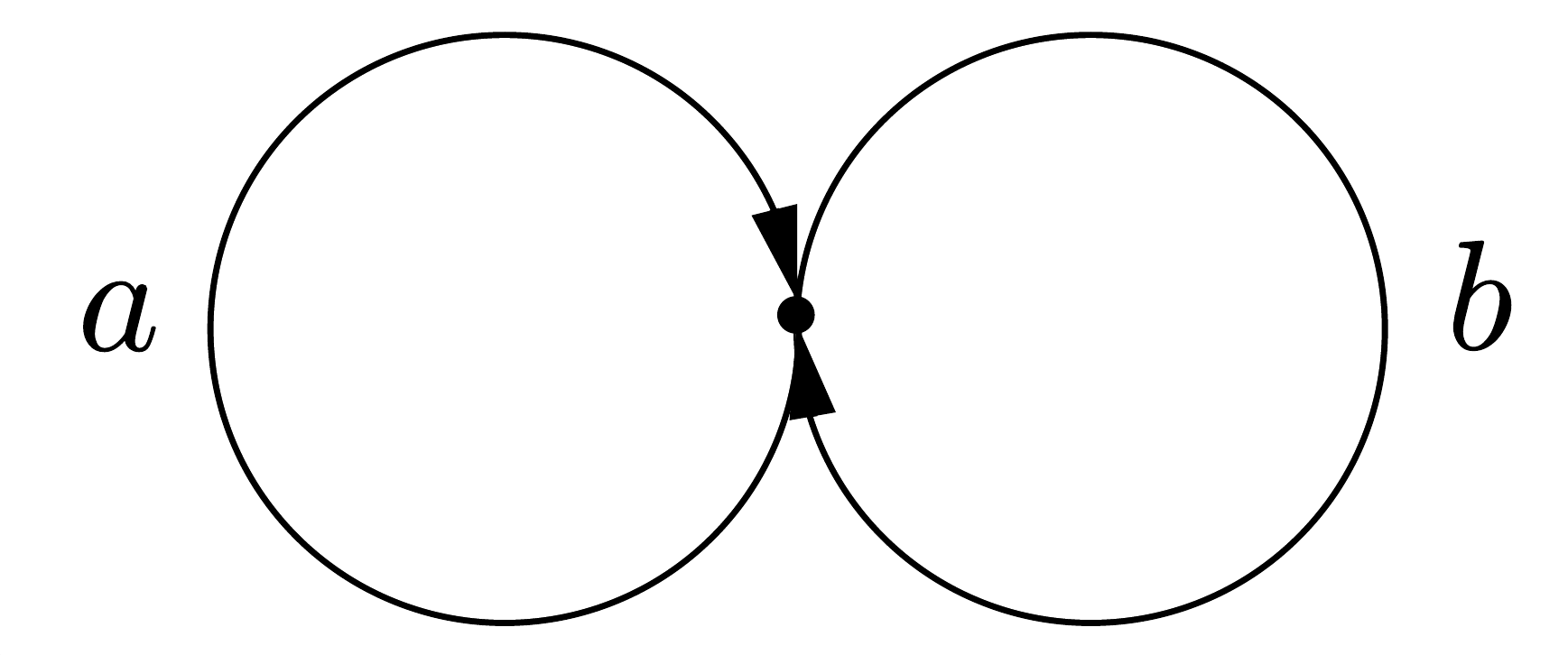}& $\langle a,b \rangle^\omega$\\
\hline
\includegraphics[trim=0 0 0 -5, width=0.4\textwidth]{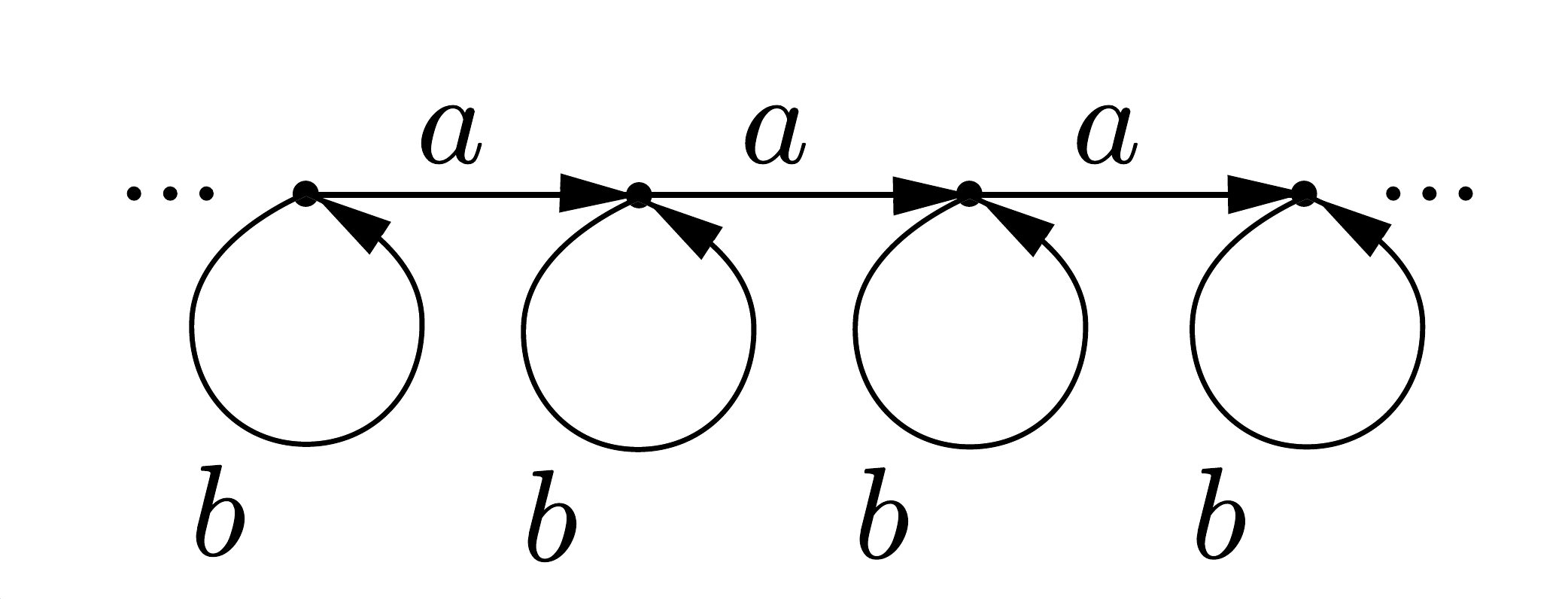}& $\langle a^nba^{-n} : n \in \mathbb Z \rangle^\omega$\\
\hline
\includegraphics[trim=0 0 0 -5, width=0.4\textwidth]{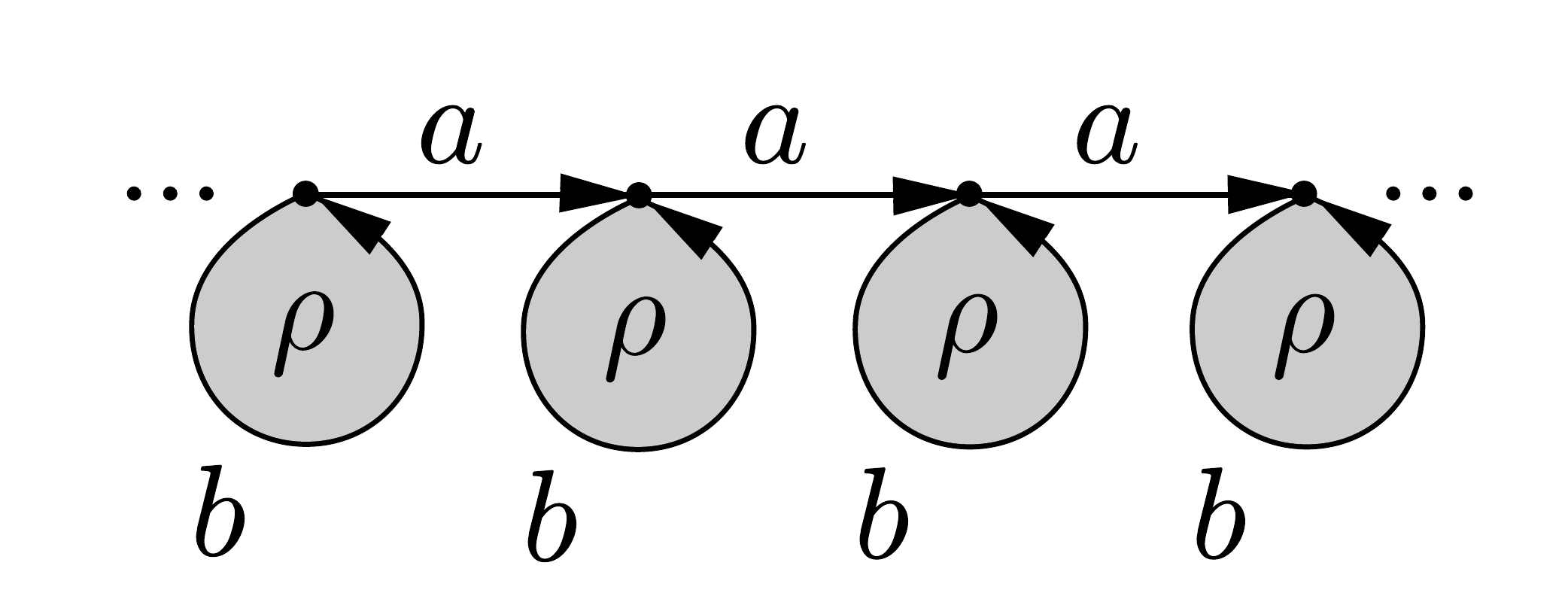}& $\langle a^n\rho a^{-n} : n \in \mathbb Z \rangle^\omega$,
the universal cover\\
\hline
\includegraphics[trim=0 0 0 -5, width=0.3\textwidth]{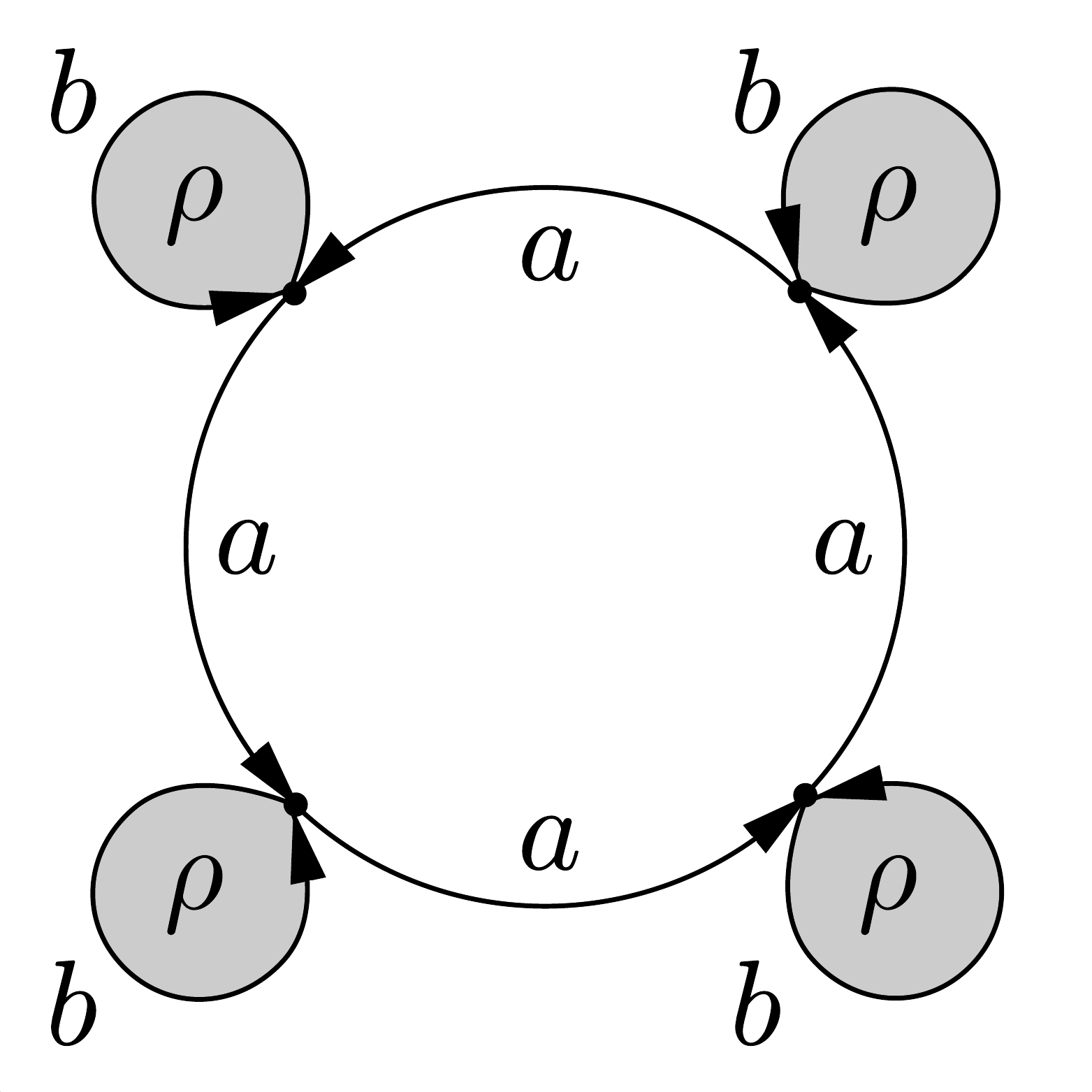}& $\langle a^k, a^n\rho a^{-n} : n \in \{1, \ldots, k\} \rangle^\omega$

$k \in \mathbb N$, $(k=4)$\\
\hline
\includegraphics[trim=0 0 0 -5, width=0.6\textwidth]{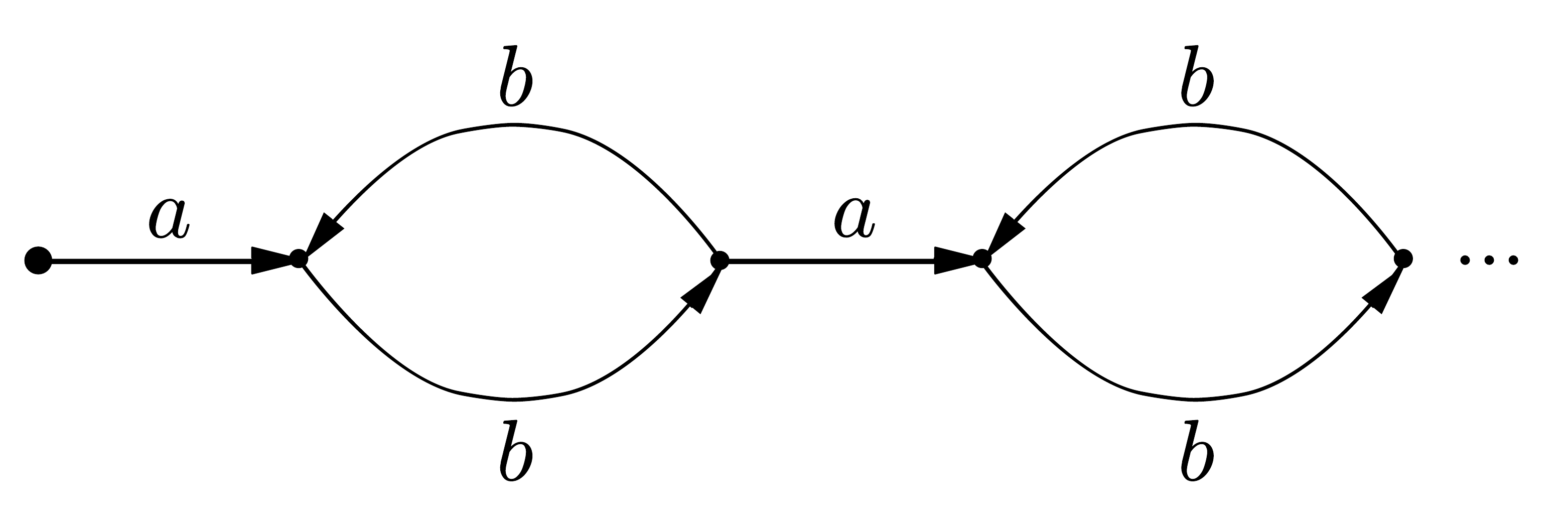}& $\langle (ab)^n ab^2 a^{-1} (ab)^{-n}: n \in \mathbb N \rangle^\omega$\\
\hline
\includegraphics[trim=0 0 0 -5, width=0.4\textwidth]{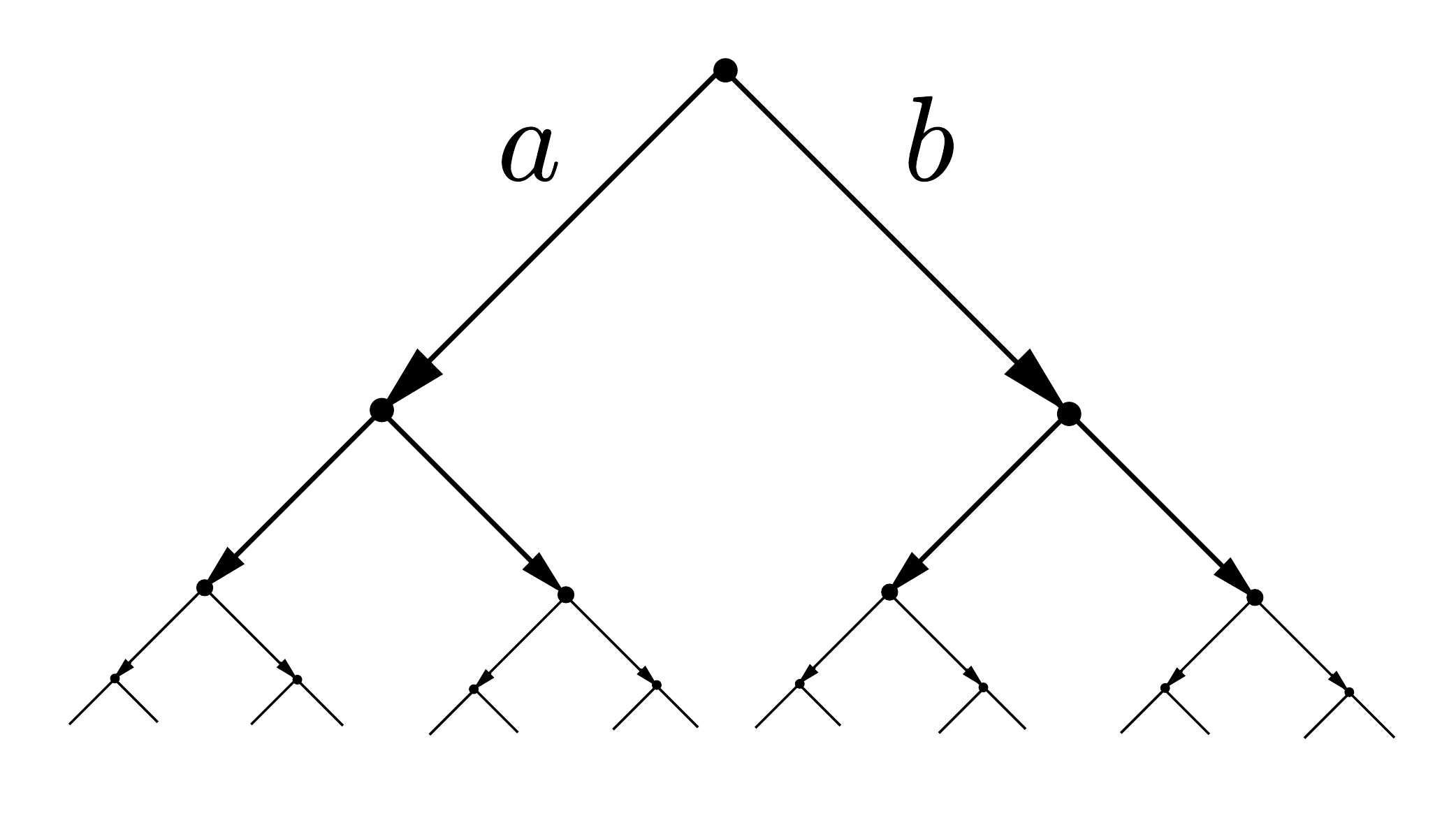}& $\{ ww^{-1} : w \in \{a,b\}^\ast \}$\\
\hline
\includegraphics[trim=0 0 0 -5, width=0.466\textwidth]{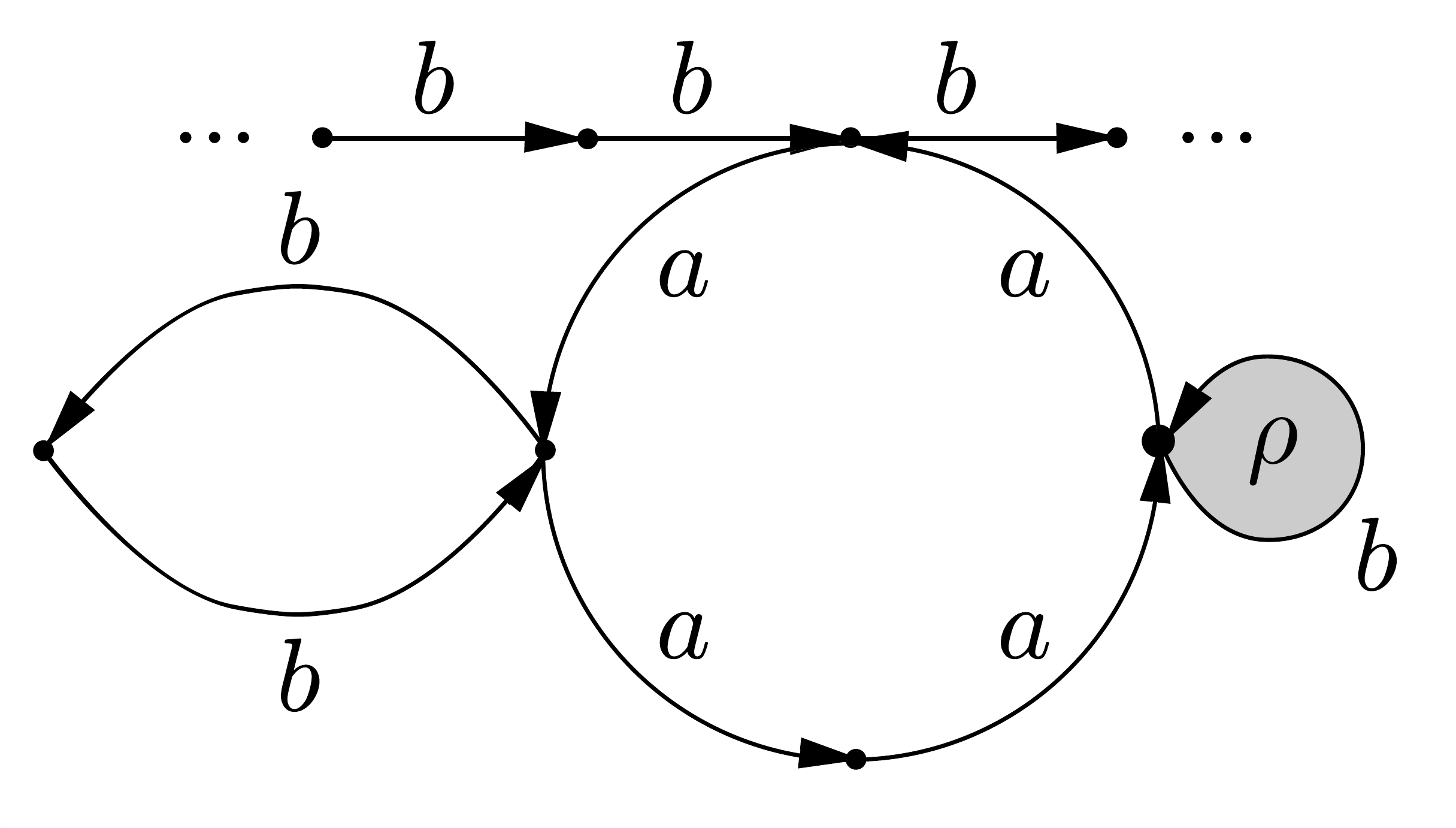}& $\langle a^4,\ \rho,\ a^2b^2a^{-2},\ ab^nb^{-n}a^{-1}: n \in \mathbb Z \rangle^\omega$\\
\hline
\includegraphics[trim=0 0 0 -5, width=0.4\textwidth]{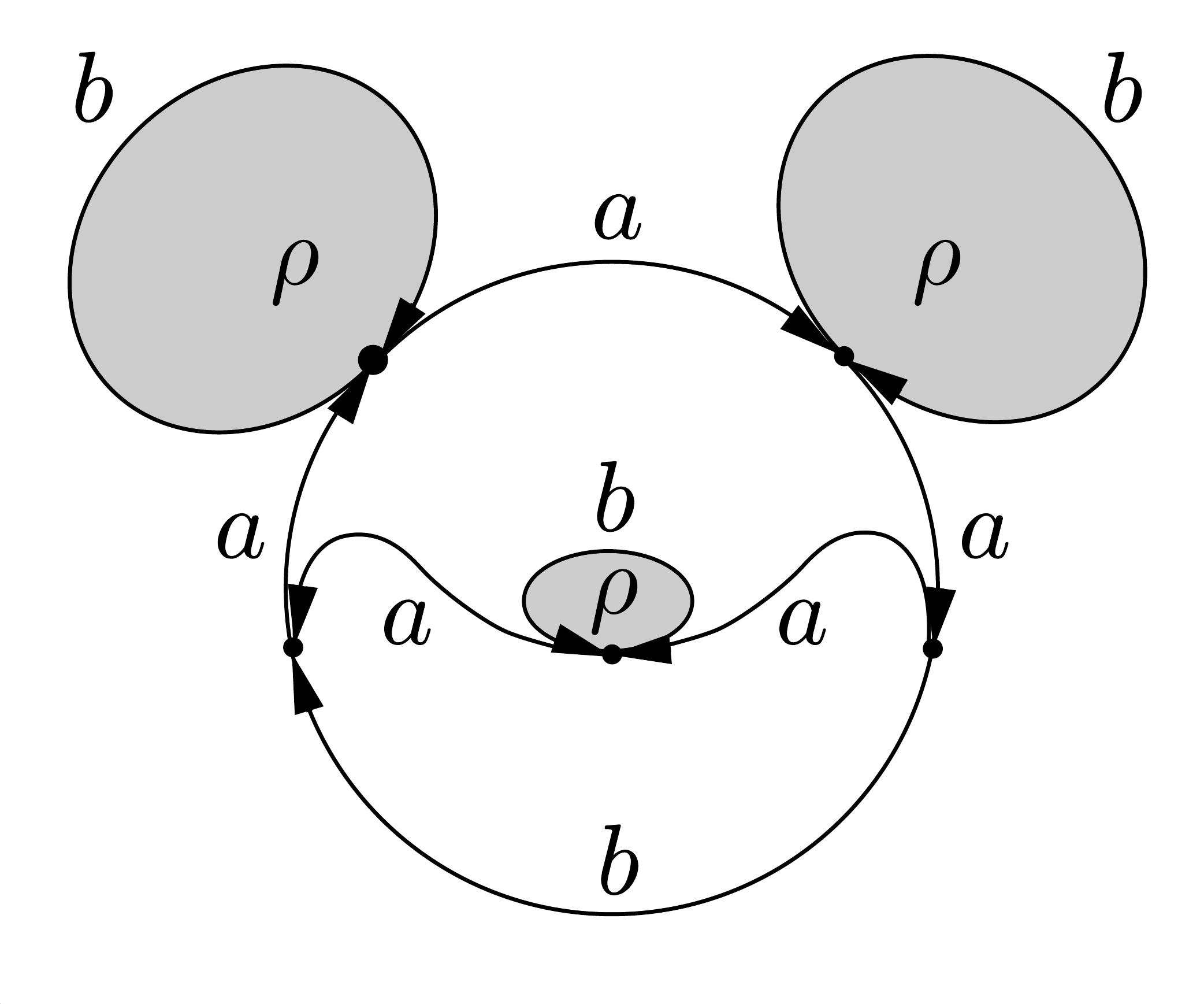}& $\langle \rho,\ a\rho a^{-1},\ a^2ba,\ a^3\rho a^2 \rangle^\omega$\\
\hline
\end{longtable}

\end{Ex}

\begin{Ex}
Regard the torus as the $2$-complex seen in Example \ref{torus}. Its loop
monoid is $\M{X,P}=\langle a,b,\rho\ |\ \rho^2=\rho, \rho = \rho
aba^{-1}b^{-1}\rangle$. We construct the unique complex
$\comp=\comp_H$ with a loop monoid $H=\langle a^{-1}b^{-1}ab,\
ab\rho a^{-1} \rangle ^\omega \leq \M{X,P}$ using the method
described in Theorem \ref{cosetalg}. 

\begin{longtable}{m{0.3\textwidth} m{0.7\textwidth}}

The flower automaton: & \includegraphics[width=0.529\textwidth]{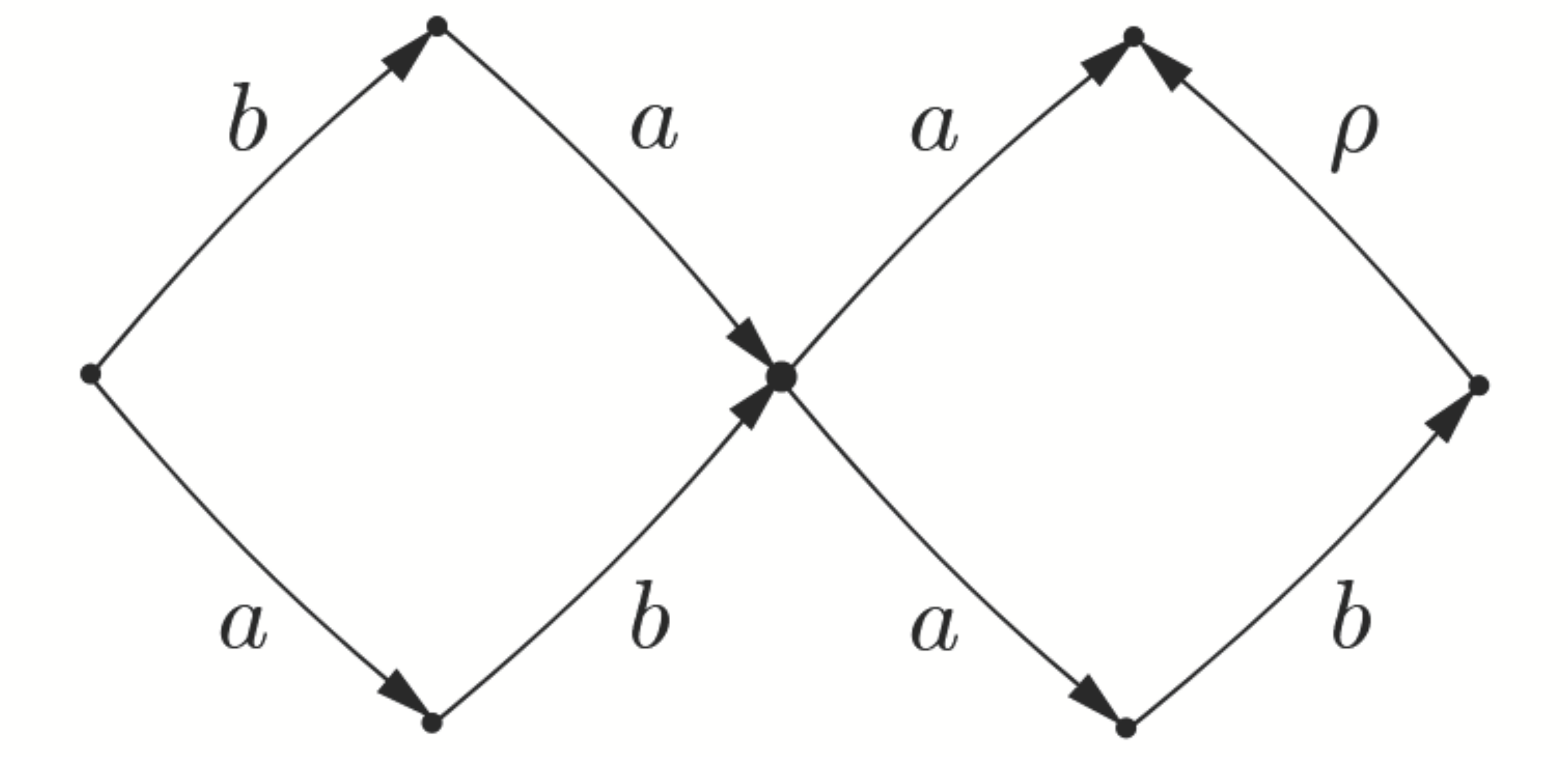}\\
Folding $a$, then expanding by
$\rho^2$: & \includegraphics[width= 0.622\textwidth]{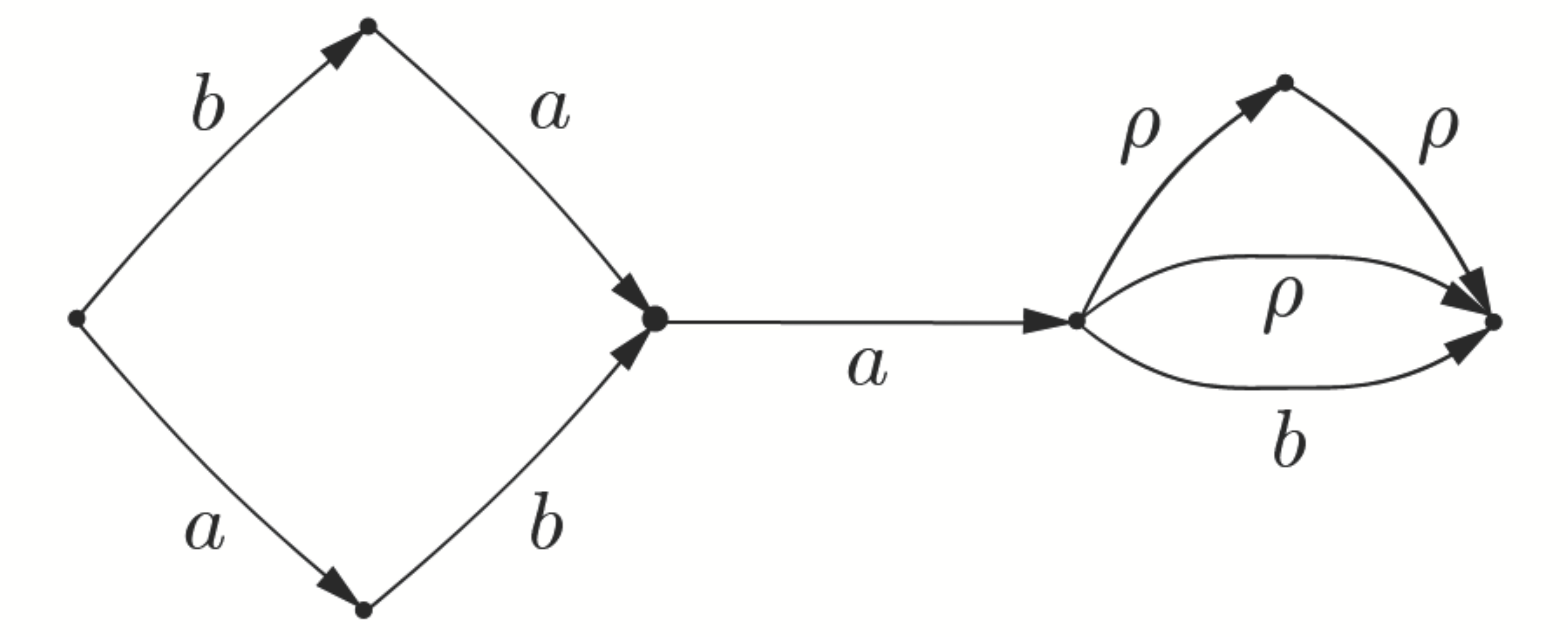}\\
Folding $\rho$, then expanding by $\rho aba^{-1}b^{-1}$, and folding
$\rho$ right away: & \includegraphics[width=0.673\textwidth]{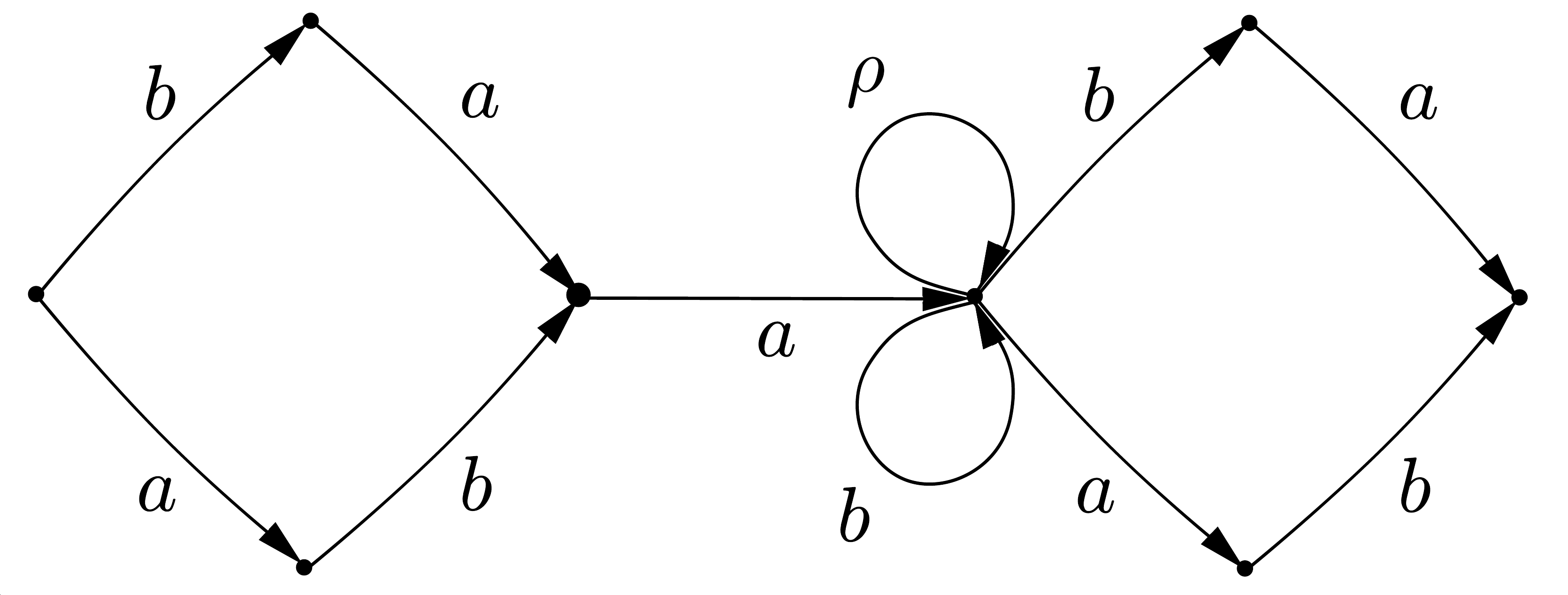}\\
Folding $b$ and then $a$, the resulting graph is complete, thus it is
$\graph_H$: & \includegraphics[width=0.594\textwidth]{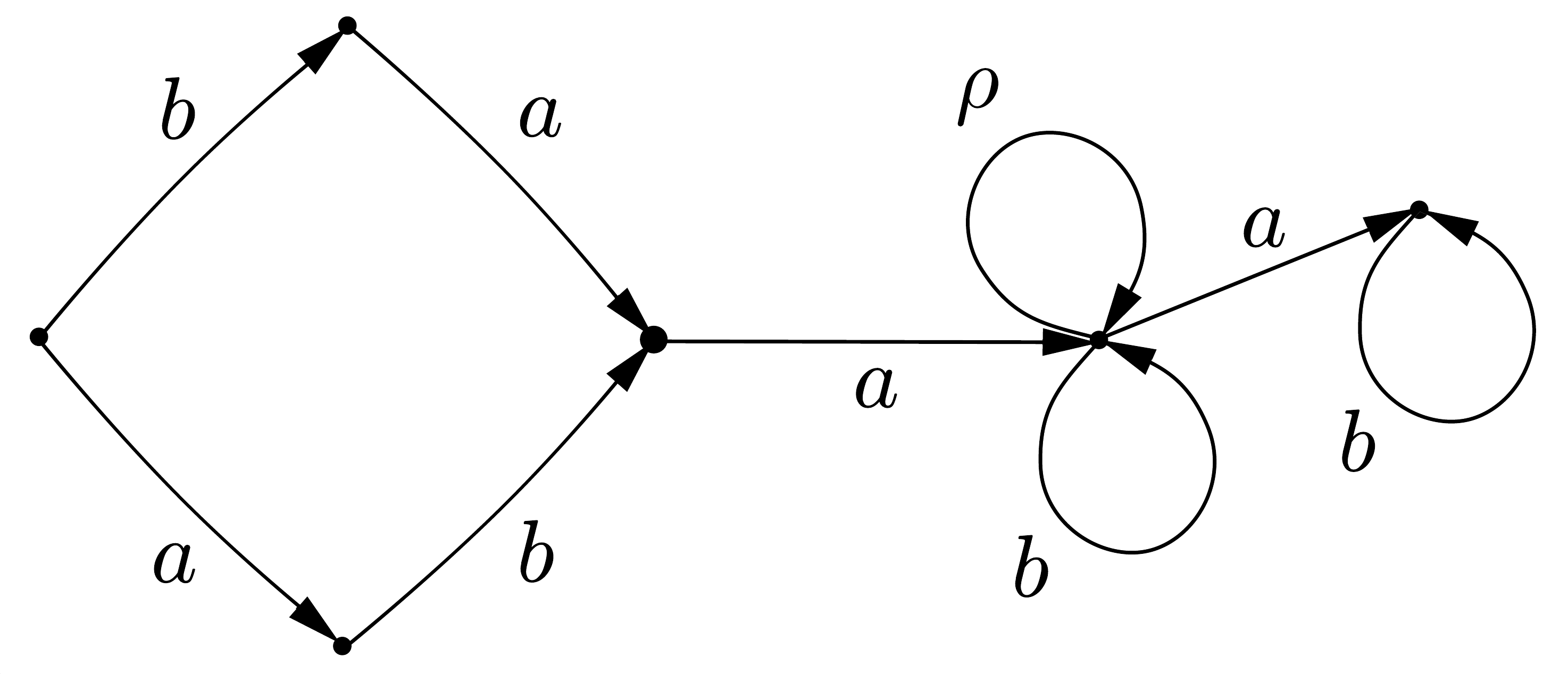}\\
The coset complex $\comp_H$: & \includegraphics[width=0.594\textwidth]{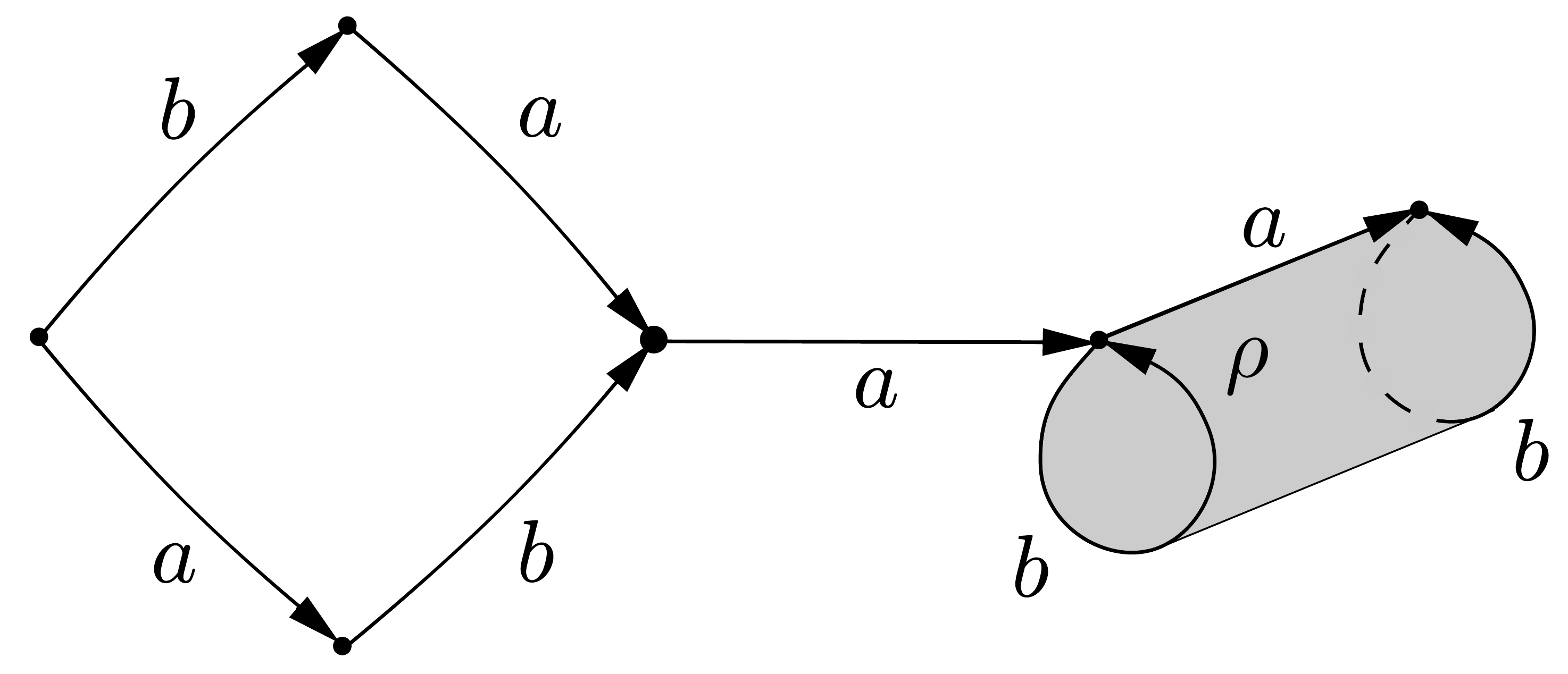}\\
\end{longtable}

\end{Ex}

\section*{Acknowledgement}

The authors are grateful to Victoria Gould for providing them with a
copy of Helen Williamson's thesis \cite{will}.

\vskip 0.5in

\end{document}